\RequirePackage{pgf,tikz,algpseudocode}
\documentclass[sn-mathphys]{sn-jnl}
\jyear{2023}%
\usepackage{bbm}
\usepackage{bm}
\usepackage{nicematrix}
\theoremstyle{thmstyleone}%
\newtheorem{theorem}{Theorem}
\newtheorem{proposition}[theorem]{Proposition}%
\newtheorem{lemma}[theorem]{Lemma}
\newtheorem{corollary}[theorem]{Corollary}
\newtheorem{assumption}[theorem]{Assumption} 

\theoremstyle{thmstyletwo}%
\newtheorem{remark}{Remark}%

\theoremstyle{thmstylethree}%

\newcommand{\wo}{\omega}
\newcommand{\LQ}[0]{\left[}
\newcommand{\RQ}[0]{\right]}

\newcommand{\setR}[0]{\mathbb{R}}
\newcommand{\PP}[0]{\mathbb{P}}
\newcommand{\E}[0]{\mathbb{E}}
\newcommand{\EAP}[0]{\widehat{\mathbb{E}}}
\newcommand{\D}[0]{\mathcal{D}}

\newcommand{\yb}{\mathbf{y}}
\newcommand{\pb}{\mathbf{p}}
\newcommand{\ybt}{\widetilde{\mathbf{y}}}
\newcommand{\pbt}{\widetilde{\mathbf{p}}}
\newcommand{\ub}{\mathbf{u}}
\newcommand{\dyb}{\mathbf{dy}}
\newcommand{\dpb}{\mathbf{dp}}
\newcommand{\dub}{\mathbf{du}}
\newcommand{\Fb}{\mathbf{F}}
\newcommand{\BB}{\mathcal{B}}
\newcommand{\unob}{\mathbf{1}}
\newcommand{\Jb}{\mathbf{J}}
\newcommand{\Fbt}{\widetilde{\mathbf{F}}}
\newcommand{\Jbt}{\widetilde{\mathbf{J}}}
\newcommand{\xb}{\mathbf{x}}
\newcommand{\jb}{\overline{j}}
\newcommand{\xbt}{\widetilde{\mathbf{x}}}
\newcommand{\dxbt}{\widetilde{\mathbf{dx}}}
\newcommand{\bb}{\mathbf{f}}
\newcommand{\vb}{\mathbf{v}}
\newcommand{\jt}{\widetilde{j}}
\newcommand{\diag}[1]{\text{diag}(#1)}
\newcommand{\Cvar}[1]{\text{CVaR}_\lambda\left( #1 \right)}
\newcommand{\Var}[1]{\text{VaR}_\lambda\left( #1 \right)}

\begin{document}

\title[A MG solver for PDE-constrained optimization with uncertain inputs]{A multigrid solver for PDE-constrained optimization with uncertain inputs}

\author[1]{Gabriele Ciaramella}\email{gabriele.ciaramella@polimi.it}

\author[2]{Fabio Nobile}\email{fabio.nobile@epfl.ch}

\author*[3]{Tommaso Vanzan}\email{tommaso.vanzan@polito.it}

\affil[1]{\orgdiv{MOX, Dipartimento di Matematica}, \orgname{Politecnico di Milano}, \country{Italy}}

\affil[2]{\orgdiv{CSQI Chair, Institute of Mathematics}, \orgname{Ecole Polytecnique F\'ed\'erale de Lausanne}, \country{Switzerland}}

\affil*[3]{\orgdiv{Dipartimento di Scienze Matematiche}, \orgname{Politecnico di Torino}, \country{Italy}}


\abstract{In this manuscript, we present a collective multigrid algorithm to solve efficiently the large saddle-point systems of equations that typically arise in PDE-constrained optimization under uncertainty, and develop a novel convergence analysis of collective smoothers and collective two-level methods. The multigrid algorithm is based on a collective smoother that at each iteration sweeps over the nodes of the computational mesh, and solves a reduced saddle-point system whose size is proportional to the number $N$ of samples used to discretized the probability space. We show that this reduced system can be solved with optimal $O(N)$ complexity. 

The multigrid method is tested both as a stationary method and as a preconditioner for GMRES on three problems: a linear-quadratic problem, possibly with a local or a boundary control, for which the multigrid method is used to solve directly the linear optimality system; a nonsmooth problem with box constraints and $L^1$-norm penalization on the control, in which the multigrid scheme is used as an inner solver within a semismooth Newton iteration; a risk-averse problem with the smoothed CVaR risk measure where the multigrid method is called within a preconditioned Newton iteration. In all cases, the multigrid algorithm exhibits excellent performances and robustness with respect to the parameters of interest.
}

\keywords{multigrid, optimization under uncertainty, random PDEs}

\pacs[MSC Classification]{65M55 - 65F10 - 65K10 - 49J55}

\maketitle

\section{Introduction}
\label{sec:intro}
In this work, we present a multigrid method to solve the saddle point system
\begin{equation}\label{eq:system_saddle_point}
\mathcal{S}\xb=\mathbf{f},
\end{equation}
where $\xb=(\yb,\ub,\pb)=(\yb_1,\dots,\yb_N,\ub,\pb_1,\dots,\pb_N)^\top$, $\mathcal{S}$ has the block structure
\begin{equation}\label{eq:matrix_saddle_point}
\mathcal{S}=\begin{pmatrix}
C_1 & & & & A_1^\top\\
& \ddots & & & &\ddots\\
& & C_N & & & & A_N^\top\\
& & & G & D_1 &\dots & D_N\\
A_1 & & & E_1\\
& \ddots & &\vdots\\
& & A_N & E_N
\end{pmatrix},
\end{equation}
and all submatrices involved represent the discretization of some differential operators. More details on each block are provided in Section \ref{Sec:quadratic}.
Matrices such as \eqref{eq:matrix_saddle_point} are often encountered while solving PDE-constrained optimization problems under uncertainty of the form
\begin{equation}\label{eq:OCP_model_introduction}
\begin{aligned}
&\min_{u\in U} \mathcal{R}\LQ Q(y(\omega),u)\RQ\\
&\text{s.t. } y(\omega)\in V\text{ satisfies}\\
& \langle e(y(\omega),u,\omega),v\rangle =0 \quad \forall v\in V,\text{ a.e. }\omega\in \Omega,
\end{aligned}
\end{equation}
where $u$ is the unknown \textit{deterministic} control, $y(\omega)$ is the state variable which satisfies a random PDE constraint expressed by $e(\cdot,\cdot,\omega)$ for almost every realization $\omega$ of the randomness, $Q$ is a real-valued quantity of interest (cost functional) and $\mathcal{R}$ is a risk measure. The vectors $\left\{\yb_j\right\}_{j=1}^N$ and $\left\{\pb_j\right\}_{j=1}^N$ are the discretizations of the state and adjoint variables $y(\omega)$ and $p(\omega)$ at the 
$N$ samples in which the random PDE constraint is collocated. The vector $\mathbf{u}$ is the discretization of the deterministic control $u$.
Problems of the form \eqref{eq:OCP_model_introduction} are increasingly employed in applications. The PDE constraints typically represent some underlying physical model whose behaviour should be optimally controlled, and the randomness in the PDE allows one to take into account the intrisinc variability or lack of knowledge on some parameters entering the model. The introduction of a risk measure in \eqref{eq:OCP_model_introduction} allows one to construct robust controls that take into account the distribution of the cost over all possible realizations of the random parameters.
Therefore, the topic has received a lot of attention in the last years, see, e.g. \cite{kouri2018optimization,Kouri_Cvar,martinez2018optimal,doi:10.1137/19M1294952,geiersbach2020stochastic,antil2021ttrisk,nobile_vanzan2,eigel2018risk,ASADPOURE20111131}.

However, few works have focused on efficient solvers for the optimality systems \eqref{eq:system_saddle_point}. A popular approach is to perform a Schur complement on $\mathbf{u}$ and solve the reduced system with a Krylov method (possibly with Conjugate Gradient), despite each iteration would then require the solution of $2N$ PDEs, with $A_j$ and $A_j^\top$ for $j=1,\dots,N$ \cite{Kourisparse}. For a full-space formulation, block diagonal preconditioners have been proposed in \cite{Kouri2018} and analyzed in \cite{nobile_vanzan}, using both an algebraic approach based on Schur complement approximations and an operator preconditioning framework.
 
In this manuscript, we design a multigrid method to solve general problems of the form \eqref{eq:system_saddle_point}, present a detailed convergence analysis which, although in a simplified setting, is nontrivial and requires technical arguments, and show how this strategy can be used for the efficient solution of three different Optimal Control Problems Under Uncertainty (OCPUU). First, we consider a linear-quadratic OCPUU and use the multigrid algorithm directly to solve the linear optimality system. Second, we consider a nonsmooth OCPUU with box constraints and $L^1$ regularization on the control. To solve such problem, we use the collective multigrid method as an inner solver within an outer semismooth Newton iteration. Incidentally, we show that the theory developed for the deterministic OCPs with $L^1$ regularization can be naturally extended to the class of OCPUU considered here.
Third, we study a risk-averse OCPUU involving the smoothed Conditional Value at Risk (CVaR) and test the performance of the multigrid scheme in the context of a nonlinear preconditioned Newton method.

The multigrid algorithm is based on a collective smoother \cite{borzi2005multigrid,borzi2009multigrid,takacs2011convergence} that, at each iteration, loops over all nodes of the computational mesh (possibly in parallel), collects all the degrees of freedom related to a node, and updates them collectively by solving a reduced saddle-point problem.
For classical (deterministic) PDE-constrained optimization problems with a distributed control, this reduced system has size $3\times 3$, thus its solution is immediate \cite{borzi2009multigrid}. In our context, the reduced problem has size $(2N+1)\times (2N+1)$, which can be large when dealing with a large number of samples. Fortunately, we show that it can be solved with optimal $O(N)$ complexity.

From the theoretical point of view, there are very few convergence analyses of collective smoothers even in the deterministic setting, namely \cite{borzi2009multigrid} based on a local Fourier analysis, and \cite{takacs2011convergence} which relies on an algebraic approach. Notably, the presence of a low-rank block matrix in the reduced optimality system (obtained by eliminating the control) as well as the need to have stiffness and mass matrices with specific structure prevents to reasonably extend the analysis of \cite{takacs2011convergence}. We therefore present in this manuscript a fully new convergence analysis of collective smoothers and two-level collective multigrid methods in a simplified setting, which also covers the deterministic setting as particular instance.

Let us remark that collective multigrid strategies have been applied to OCPUU in \cite{Borzi,Borzi2} and in \cite{rosseel2012optimal}. This manuscript differs from the mentioned works since, on the one hand, \cite{Borzi,Borzi2} considers a \textit{stochastic} control $u$, therefore for (almost) every realization of the random parameters a different control $u(\omega)$ is computed through the solution of a standard deterministic OCP. On the other hand, \cite{rosseel2012optimal} considers a stochastic Galerkin discretization, and hence the correspoding optimality system has a structure which is very different from \eqref{eq:matrix_saddle_point}.

The multigrid algorithm presented here assumes that all state and adjoint variables are discretized on the same finite element mesh. The control can instead live on a subregion of the computational mesh, so that the algorithm is applicable also to optimization problems with local or boundary controls.

Finally, we remark that the multigrid solver proposed is based on a hierarchy of spatial discretizations corresponding to different levels of approximation, but the discretization of the probability space remains fixed, that is, the number of samples remains constant across the multigrid hierarchy. The extension of the multigrid algorithm to coarsening procedures also in the probability space will be the subject of future endeavours. We hint at possible approaches and challenges in Section \ref{Sec:collective} (see Remark \ref{remark:multilevel}). Nevertheless, we stress that the multigrid algorithm can already be incorporated within outer optimization routines that take advantange of different levels of approximations of the probability space, see, e.g., \cite{Kourisparse,kouri2014multilevel,nobile_vanzan2}.

The rest of the manuscript is organized as follows. In Section \ref{Sec:quadratic} we introduce the notation, a classical linear-quadratic OCPUU, and interpret \eqref{eq:matrix_saddle_point} as the matrix associated to the optimality system of a discretized OCPUU. Section \ref{Sec:collective} presents the collective multigrid algorithm, discusses implementation details and develops the convergence analysis. Further, the algorthm is numerically tested on the linear-quadratic OCPUU. In Section \ref{Sec:l1}, we consider a nonsmooth OCPUU with box constraints and a $L^1$ regularization on the control.
Section \ref{Sec:risk} deals with a risk-averse OCPUU.
For each of these cases, we first show how the multigrid approach can be integrated into the solution process, by detailing concrete algorithms, and then we present extensive numerical experiments to show the efficiency of the proposed framework.
Finally, we draw our conclusions in Section \ref{Sec:conc}.

\section{A linear-quadratic optimal control problem under uncertainty}\label{Sec:quadratic}
Let $\D\subset \mathbb{R}^d$ be a Lipschitz bounded domain, $V\subset L^2(\D)$ a Sobolev space (e.g. $H^1(\D)$ equipped with suitable boundary conditions), and $(\Omega,\mathcal{F},\PP)$ a complete probability space. Given a function $u$ belonging to a Hilbert space $U$, we consider the linear elliptic random PDE
\begin{equation}\label{eq:random_PDE_weak}
a_\omega(y,v)=\langle \BB u,v\rangle,\forall v\in V,\quad \PP\text{-a-e. } \omega\in \Omega,
\end{equation}
where $a_{\omega}(\cdot,\cdot):V\times V\rightarrow \mathbb{R}$ is a bilinear form and $\langle\cdot,\cdot \rangle$ denotes the duality between $V$ and $V^\prime$. $\BB:U\rightarrow V^\prime$ is a continuous control operator allowing possibly for a local control (i.e. a control acting only on a subset $\D_0\subset \D)$ or a boundary control (i.e. a control acting as Neumann condition on a subset of $\partial \D$).
To assure uniqueness and sufficient integrability of the solution of \eqref{eq:random_PDE_weak}, we make the following additional assumption. 
\begin{assumption}\label{ass:diff}
There exist two random variables $a_{\min}(\omega)$ and $a_{\max}(\omega)$ such that
\[0<a_{\min}(\omega)\|v\|^2_V\leq a_\omega(v,v)\leq a_{\max}(\omega)\|v\|^2_V,\quad \forall v\in V,\ \PP\text{-a.e. }\omega \in \Omega,\]
and further $a^{-1}_{\min}$ and $a_{\max}$ are in $L^p(\Omega)$ for some $p\geq 4$.
\end{assumption}
Under Assumption \ref{ass:diff}, it is well-known (see, e.g., \cite{lord_powell_shardlow_2014,Scheichl}) that \eqref{eq:random_PDE_weak} admits a solution in $V$ for $\PP\text{-a.e. } \omega$, and the solution $y$, interpreted as a $V$-valued random variable $y:\omega\in \Omega\mapsto y(\omega)\in V$, lies in the Bochner space $L^q(\Omega;V)$, $q\leq p$, \cite{cohn2013measure}. 
We often use the shorthand notation $y_\omega=y(\cdot,\omega)$ when the dependence on $x$ is not needed, or $y_{\omega}(u)$ if we wish to highlight the dependence on the control function $u$.

In this manuscript, we consider the minimization of functionals constrained by \eqref{eq:random_PDE_weak}.
Let us first focus on the linear-quadratic problem
\begin{equation}\label{eq:quadratic_OCP}
\begin{aligned}
&\min_{u\in U,y\in L^2(\Omega;V)} \frac{1}{2}\E\LQ \|\mathcal{I}y_\omega-y_d\|^2_{L^2(\D)}\RQ +\frac{\nu}{2}\|u\|^2_{U},\\
&\quad\text{subject to}\\
&a_\omega(y_\omega,v)=\langle \BB u+f,v\rangle,\quad \forall v \in V,\ \PP\text{-a.e. } \omega\in \Omega,
\end{aligned}
\end{equation}
where $y_d\in L^2(\D)$ is a target state, $f\in V^\prime$, $\E:L^1(\Omega)\rightarrow \setR$ is the expectation operator, $\nu >0$, and $\mathcal{I}$ is the embedding operator from $V$ to $L^2(\D)$.

Introducing the linear control-to-state map $S: g\in V^\prime \rightarrow y_\omega(g)\in L^2(\Omega;V)$, the reduced formulation of \eqref{eq:quadratic_OCP} is
\begin{equation}\label{eq:quadratic_OCP_reduced}
\min_{u\in U} \frac{1}{2}\E \LQ \|\mathcal{I}S(\mathcal{B} u+f)-y_d\|^2_{L^2(\D)}\RQ +\frac{\nu}{2}\|u\|^2_{U}.
\end{equation}
Existence and uniqueness of the minimizer of \eqref{eq:quadratic_OCP_reduced} follows directly from standard variational arguments \cite{lions1971optimal,hinze2008optimization,troltzsch2010optimal,kouri2018optimization}.
Furthermore, due to Assumption \ref{ass:diff}, the optimal control $\overline{u}$ satisfies the variational equality
\begin{equation}\label{eq:optimality_condition_quadratic_exact}
(\nu \overline{u} -\Lambda_U \BB^\star S^\star \mathcal{I}^\star(y_d-S(\mathcal{B}\overline{u}+f)),v)_{U}=0,\quad \forall v \in U,
\end{equation}
where $\Lambda_U$ is the Riesz operator of $U$. The adjoint operator $S^\star: L^2(\Omega;V^\prime)\rightarrow  V$ 
is characterized by $S^\star z=\E\LQ p\RQ$ where $p=p_{\omega}(x)$ is the solution of the adjoint equation
\begin{equation}\label{eq:adjoint_equation}
a_\omega(v,p_\omega)=\langle z(\omega),v\rangle,\quad \forall v \in V,\ \PP\text{-a-e. } \omega\in \Omega.
\end{equation}
The optimality condition \eqref{eq:optimality_condition_quadratic_exact} can thus be formulated as the optimality system
\begin{equation}\label{eq:optimality_system}
\begin{aligned}
& a_\omega(y_\omega,v)=\langle \BB \overline{u}+f,v\rangle,\quad \forall v\in V,\quad \PP\text{-a-e. } \omega\in \Omega,\\
& a_\omega(v,p_\omega)=\langle \mathcal{I}^\star(y_d-y_\omega),v\rangle,\quad \forall v \in V,\ \PP\text{-a-e. } \omega\in \Omega,\\
& (\nu \overline{u}- \Lambda_U \BB^\star\E\LQ p_\omega\RQ,v)_{U}=0,\quad \forall v \in U.
\end{aligned}
\end{equation}
To solve numerically \eqref{eq:quadratic_OCP}, we replace the exact expectation operator $\E$ of the objective functional by a quadrature formula $\EAP$ with $N$ nodes $\left\{\omega_i\right\}_{i=1}^N$ and positive weights $\left\{\zeta_i\right\}_{i=1}^N$, namely
\[ \E\LQ X\RQ \approx \EAP\LQ X\RQ:= \sum_{i=1}^N \zeta_i X(\omega_i),\quad \text{with}\quad \sum_{i=1}^N \zeta_i=1.\]
Common quadrature formulae are Monte Carlo, Quasi-Monte Carlo and Gaussian formulae. The latter requires that the probability space can be parametrized by a (finite or countable) sequence of random variables $\left\{\chi_j\right\}_j$, each with distribution $\mu_j$, and the existence of a complete basis of tensorized $L^2_{\mu_j}$-orthonormal polynomials.
Hence for the semi-discrete OCP, the $\PP$-a.e. PDE-constraint is naturally collocated onto the nodes of the quadrature formula.

Concerning the space domain, we consider a family of regular triangulations $\left\{\mathcal{T}_h\right\}_{h>0}$ of $\D$, and a Galerkin projection onto a conforming finite element space $V^h\subset V$ of continuous piecewise polynomial functions of degree $r$ over $\mathcal{T}_h$.
$N_h$ is the dimension of $V^h$ and $\left\{\phi_i\right\}_{i=1}^{N_h}$ is a nodal Lagrangian basis. We discretize the state and adjoint variables on the same finite element space.
The control variable is discretized on the finite element space $U_h=\text{span}\left\{\psi_i\right\}_{i=1}^{N_u}$, where $N_u$ is possibly strictly smaller than $N_h$ in case of a local or a boundary control.

Once fully discretized, \eqref{eq:optimality_system} can be expressed as
\begin{equation}\label{eq:full_space_optimality_system}
\begin{pmatrix}
M & & & & A_1^\top\\
& \ddots & & & &\ddots\\
& & M & & & & A_N^\top\\
& & & \nu M_U & -\zeta_1 B^\top &\dots & -\zeta_N B^\top\\
A_1 & & & -B\\
& \ddots & &\vdots\\
& & A_N & -B
\end{pmatrix}
\begin{pmatrix}
\yb_1\\ \vdots \\\yb_N\\ \mathbf{u}\\
\pb_1\\ \vdots \\ \pb_N
\end{pmatrix}=
\begin{pmatrix}
M \yb_d\\ \vdots \\
M \yb_d\\ \mathbf{0} \\ M\mathbf{f}\\ \vdots \\ M\mathbf{f}
\end{pmatrix},
\end{equation}
where $A_j$ are the stiffness matrices associated to the bilinear forms $a_{\omega_j}(\cdot,\cdot)$, $M$ and $M_U$ are mass matrices corresponding to the finite element spaces $V_h$ and $U_h$, $B$ is the discretization of the control operator, $\yb_d$ and $\mathbf{f}$ are the finite element discretizations of $y_d$ and $f$ respectively, while $\yb_j$ and $\pb_j$ are the discretizations of $y_{\omega_j}$ and $p_{\omega_j}$. 
Notice that the matrix in \eqref{eq:full_space_optimality_system} could be symmetrized by multipling the first and the last $N$ rows by the weights $\left\{\zeta_i\right\}_{i=1}^N$. This would be also consistent with the theoretical interpretation of the blocks of the saddle point system as discretizations of continuous inner products. From the numerical point of view, we have not observed relevant advantanges in maintaining the weights. Since for more general problems (see, e.g., Sec. \ref{Sec:l1}) the symmetry of the saddle point system cannot be recovered by multiplying some equations by the quadrature weights, we do not consider the symmetrized version in this work.

\section{Collective multigrid scheme}\label{Sec:collective}
In this section, we describe the multigrid algorithm to solve the full space optimality system \eqref{eq:full_space_optimality_system}. First, we consider a distributed control, so that $u$ lives on the whole computational mesh and $B=M$. Local and boundary controls are discussed at the end of the section.
Second, for the sake of generality, we consider the more general matrix \eqref{eq:matrix_saddle_point}, so that our discussion covers also the different saddle-point matrices obtained in Sections \ref{Sec:l1} and \ref{Sec:risk}. 

For each node of the triangulation, let us introduce the vectors $\ybt_i$ and $\pbt_i$,
\[\ybt_i=\begin{pmatrix}
(\yb_1)_i\\  \vdots\\(\yb_N)_i
\end{pmatrix}\in \setR^{N},
\quad \pbt_i=\begin{pmatrix}
(\pb_1)_i\\\vdots\\(\pb_N)_i
\end{pmatrix} \in \setR^{N}, \quad i=1,\dots,N_h,
\]
which collect the degrees of freedom associated to the $i$-th node, the scalar $u_i=(\mathbf{u})_i$,
and the restriction operators $R_i \in \mathbb{R}^{(2N+1) \times  ((2N+1) N_h)}$ such that 
\begin{equation}\label{eq:restriction_action}
R_i \begin{pmatrix}
\yb \\ \mathbf{u}\\\pb 
\end{pmatrix}=\begin{pmatrix}
\ybt_i\\ u_i\\\pbt_i
\end{pmatrix}=:\xb_i.
\end{equation}
The prolongation operators are $P_i:=R_i^\top$, while the reduced matrices $\widetilde{S}_i:=R_iSP_i\in \setR^{(2N+1)\times (2N+1)}$ represent a condensed saddle-point matrix on the $i$-th node, and satisfy

\begin{equation}\label{eq:reduced_matrices}
\widetilde{S}_i=\begin{pmatrix}
\diag{\mathbf{c}_i} & 0 &\diag{\mathbf{a}_i}\\
0 & (G)_{i,i} & \mathbf{d}_i^\top\\
\diag{\mathbf{a}_i} & \mathbf{e}_i & 0
\end{pmatrix}
\end{equation}
with $\mathbf{c}_i:=((C_1)_{i,i},\dots,(C_N)_{i,i})^\top$, $\mathbf{a}_i:=((A_1)_{i,i},\dots,(A_N)_{i,i})^\top$, $\mathbf{e}_i=((E_1)_{i,i},\dots,(E_N)_{i,i})^\top$, $\mathbf{d}_i=((D_1)_{i,i},\dots,(D_N)_{i,i})^\top$, where $\diag{\mathbf{v}}$ denotes a diagonal matrix with the components of $\mathbf{v}$ on the main diagonal.\\
Given an initial vector $\xb^0$, a Jacobi-type collective smoothing iteration computes for $n=1,\dots,n_1$,
\begin{equation}\label{eq:collective_smoothing}
\xb^n=\xb^{n-1}+ \theta\sum_{i=1}^{N_h} P_i \widetilde{S}_i^{-1}R_i\left(\bb-S\xb^{n-1}\right),
\end{equation}
where $\theta \in (0,1]$ is a damping parameter. Gauss-Seidel variants can straightforwardly be defined.
Next, we consider a sequence of meshes $\left\{\mathcal{T}_{h_\ell}\right\}_{\ell=\ell_{\min}}^{\ell_{\max}}$, which we assume for simplicity to be nested, and restriction and prolongator operators $R_{\ell-1}^\ell$, $P_{\ell-1}^{\ell}$ which map between grids $\mathcal{T}_{h_{\ell-1}}$ and $\mathcal{T}_{h_{\ell}}$. In the numerical experiments, the coarse matrices are defined recursively in a Galerkin fashion starting from the finest one, namely $S_\ell:=R^{\ell+1}_{\ell}S_{\ell+1}P^{\ell+1}_{\ell}$ for $\ell\in\left\{1,\dots, \ell_{\max}-1\right\}$. Nevertheless it is obviously possible to define $S_\ell$ as the discretization of the continuous saddle-point system onto the mesh $\mathcal{T}_{h_\ell}$.
With this notation, the V-cycle collective multigrid is described by Algorithm \ref{Alg:algorithm-Vcycle}, which can be repeated until a certain stopping criterium is satisfied. We used the notation \textit{Collective\_Smoothing$(\cdot,\cdot,\cdot)$} to denote possible variants of \eqref{eq:collective_smoothing} (e.g. Gauss-Seidel).
\begin{algorithm}[t]
\setlength{\columnwidth}{\linewidth}
\caption{V-cycle Collective Multigrid Algorithm - V-cycle($\mathbf{x}^{0}$,$\mathbf{f}$,$\ell$)}\label{Alg:algorithm-Vcycle}
\begin{algorithmic}[1]
\If{$\ell=\ell_{\min}$,}
\State set $\mathbf{x}^{0}=S^{-1}_{\ell_{\min}}\mathbf{f}.$  $\qquad \qquad \quad \quad \; \:$ $\;\;\;\;\;\;\;\; \; \; \; $ (direct solver)
\Else
\State $\xb^{n_1}$=Collective\_Smoothing$(\xb^0,S_\ell,n_1)$ $\;$ ($n_1$ steps of coll. smoothing)
\State $\mathbf{r} = \mathbf{f}-S_{\ell}\xb^{n_1}$ $\qquad \qquad \qquad \quad \; \,$ $\;\;\;\;\;\;\;\; \; \; \; \;$(compute the residual)
\State ${\bf e}_c=$V-cycle(${\bf 0}, R^\ell_{\ell-1}{\bf r},\ell-1$).$\quad \; \; \:$ $\;\;\;\;\;\; \; \; \; \;$(recursive call)
\State $\xb^{0} = \xb^{n_1}+P^\ell_{\ell-1}{\bf e}_c$ $\qquad \qquad \quad \; \: $ $\;\;\;\;\;\;\; \; \; \; \;$(coarse correction)
\State $\xb^{n_2}$=Collective\_Smoothing$(\xb^0,S_\ell,n_2)$ $\;$ ($n_2$ steps of coll. smoothing)
\State Set $\xb^{0}=\xb^{n_2}$ $\qquad \qquad \qquad \quad \; \; \; \; \; \; \; \; \; \; \; \; \; \; \;$ (update)
\EndIf\\
\Return $\xb^{0}$.
\end{algorithmic}
\end{algorithm}

Notice that \eqref{eq:collective_smoothing} requires to invert the matrices $S_i$ for each computational node. We now show that this can be done with optimal $O(N)$ complexity. Indeed, performing a Schur complement on $u_i$, the system $\widetilde{S}_i\xb_i=\bb_i$, with $\bb_i=(\bb_{p_i},b_{u_i},\bb_{y_i})^\top$ can be solved exclusively computing inverses of diagonal matrices and scalar products between vectors through
\begin{equation}\label{eq:solution_reduced_system}
\begin{aligned}
u_i&=\frac{b_{u_i}+\mathbf{d}_i^\top(\diag{\mathbf{a}_i}^{-1}\diag{\mathbf{c}_i}\diag{\mathbf{a}_i}^{-1}\bb_{y_i}-\diag{\mathbf{a}_i}^{-1}\bb_{p_i})}{(G)_{i,i}+\mathbf{d}_i^\top \diag{\mathbf{a}_i}^{-1}\diag{\mathbf{c}_i}\diag{\mathbf{a}_i}^{-1}\mathbf{e}_i},\\
\ybt_i&=(\diag{\mathbf{a}_i})^{-1}(\bb_{y_i}-\mathbf{e}_i u_{i}),\\
\pbt_i&=(\diag{\mathbf{a}_i})^{-1}(\bb_{p_i}-\diag{\mathbf{c}_i} \ybt_i).
\end{aligned}
\end{equation}
Notice that we should guaranteee that  $\diag{\mathbf{a}_i}$ admits an inverse and that $(G)_{i,i}+\mathbf{d}_i^\top \diag{\mathbf{a}_i}^{-1}\diag{\mathbf{c}_i}\diag{\mathbf{a}_i}^{-1}\mathbf{e}_i\neq 0$. This has to be verified case by case, so we now focus on the specific matrix \eqref{eq:full_space_optimality_system}. On the one hand, the vectors $\mathbf{a}_i$ are strictly positive componentwise, since $(\mathbf{a}_i)_j=a_{\omega_j}(\phi_i,\phi_i)>0$ $\forall i=1,\dots,N_h$, $j=1,\dots,N$ (due to Assumption \ref{ass:diff}). On the other hand, $(G)_{i,i}=\int_\D\psi^2_i(x)\ dx >0$, while
a direct calculation shows that  
\[\mathbf{d}_i^\top \diag{\mathbf{a}_i}^{-1}\diag{\mathbf{c}_i}\diag{\mathbf{a}_i}^{-1}\mathbf{e}_i=(M)^3_{i,i} \sum_{j=1}^N \zeta_j  (A_j)^{-2}_{i,i}>0,\]
which implies that the denominator in the first equation of \eqref{eq:solution_reduced_system} is strictly positive.

The collective smoother can be easily adjusted to accomodate local or boundary controls as discussed in \cite{borzi2003multigrid} for deterministic OCPs.  
For all nodes $i$ for which a control basis function is present, the smoothing procedure remains that of \eqref{eq:solution_reduced_system}. For all others computational nodes for which there is not a control basis function associated, the smoothing procedure becomes
\begin{equation}\label{eq:solution_reduced_system2}
\begin{aligned}
\ybt_i&=(\diag{\mathbf{a}_i})^{-1}\bb_{yi},\\
\pbt_i&=(\diag{\mathbf{a}_i})^{-1}(\bb_{pi}-\diag{\mathbf{c}_i} \ybt_i),
\end{aligned}
\end{equation}
which is consistently obtained from \eqref{eq:solution_reduced_system} setting $u_i=0$.

To conclude this section, we remark that the computational complexity of the smoothing procedure is of order $O(N_h N)$, thus linear with respect to the size of the saddle point-system. Provided that the V-cycle algorithm requires a constant number of iterations to converge as the number of levels increases, and that $N$ is not too large (so that the cost of the coarse solver is not dominant), the complexity of the multigrid algorithm can also be considered linear. In the next numerical experiments sections, we show indeed that the number of iterations remains constant for several test cases.
\begin{remark}[Extension to a hierarchy of samples]\label{remark:multilevel}
The multigrid algorithm presented is based on a hierarchy of spatial discretizations. However, the sample to discretize the probability space remains fixed among the levels. If one relies on the stochastic collocation method to discretize the probability space, it is possible to envisage a multigrid algorithm that also involves a coarsening of the sample size, since for each sample set one could consider the associated stable interpolator which can then be evaluated onto a coarser or finer set of samples. Nevertheless, it is not clear at the moment the interplay between the smoothing and coarsening procedures, which is key for the efficient behaviour of a multigrid scheme.
Future endeavours will investigate this interesting direction. For the rest of the manuscript we restrict oursevels to a hierarchy of spatial discretizations since on the one hand, the multigrid algorithm can already be embedded in other outer optimization algorithms that involve a hierarchy of samples \cite{Kourisparse,kouri2014multilevel,nobile_vanzan2,van2019robust}. On the other hand, the reduced system can be solved with optimal $O(N)$ linear complexity, so that a coarsening in the number of samples may be superfluous.
\end{remark}

\subsection{Convergence analysis}
In this subsection, we present a convergence analysis of the collective multigrid algorithm in a simplified setting. Let $\mathcal{D}=(0,1)$, and consider the random PDE
\begin{equation}\label{eq:simplified}
\eta(\wo)\int_0^1 \partial_x y(x,\wo) \partial_x v(x)\;dx = \int_0^1 (f(x)+u(x))v(x)\; dx,\forall v\in V,\; \mathbb{P}\text{-a.e.}\; \wo \in \Omega,
\end{equation}
where $\eta:\Omega\rightarrow \mathbb{R}^+$ is a positive valued random variable such that $\mathbb{E}\LQ \eta^{-2}\RQ <\infty$. Our goal is to minimize the objective functional of \eqref{eq:quadratic_OCP} constrained by \eqref{eq:simplified}. A discretization using finite differences and with $N$ Monte Carlo samples leads to the optimality system
\begin{equation}\label{eq:full_space_optimality_system_simplified}
\begin{pmatrix}
\medmuskip=-1mu
\thinmuskip=-1mu
\thickmuskip=-1mu
\nulldelimiterspace=0.9pt
\scriptspace=0.9pt 
\arraycolsep0.9em 
\frac{\widetilde{I}}{N} & & & & \frac{\eta_1(\wo)}{N} A\\
& \ddots & & & &\ddots\\
& & \frac{\widetilde{I}}{N} & & & & \frac{\eta_N(\wo)}{N} A\\
& & & \nu \widetilde{I} & -\frac{\widetilde{I}}{N} &\dots & -\frac{\widetilde{I}}{N}\\
\frac{\eta_1(\wo)}{N} A & & & -\frac{\widetilde{I}}{N}\\
& \ddots & &\vdots\\
& & \frac{\eta_N(\wo)}{N} A\ & -\frac{\widetilde{I}}{N}
\end{pmatrix}
\begin{pmatrix}
\yb_1\\ \vdots \\\yb_N\\ \mathbf{u}\\
\pb_1\\ \vdots \\ \pb_N
\end{pmatrix}=
\begin{pmatrix} \frac{\yb_d}{N}\\ \vdots \\
 \frac{\yb_d}{N}\\ \mathbf{0} \\ \frac{\mathbf{f}}{N}\\ \vdots \\ \frac{\mathbf{f}}{N}
\end{pmatrix},
\end{equation}
where $A$ is the tridiagonal matrix associated with the 1D Laplacian, with $2/h^2$ on the main diagonal, and $-1/h^2$ on the two adjacent diagonals, $h$ being the mesh size, $\widetilde{I}\in \mathbb{R}^{N_h\times N_h}$ is the identity matrix, and, compared to \eqref{eq:full_space_optimality_system}, the first and last blocks of $N$ equations are divided by $\frac{1}{N}$ to get a symmetric system. To perform our analysis, we first eliminate the variable $\mathbf{u}$, and obtain the reduced matrix
\begin{equation}\label{eq:system_reduced}
\begin{pmatrix}
\medmuskip=-1mu
\thinmuskip=-1mu
\thickmuskip=-1mu
\nulldelimiterspace=0.9pt
\scriptspace=0.9pt 
\arraycolsep0.9em 
\frac{\widetilde{I}}{N} & & & & \frac{\eta_1(\wo)}{N} A\\
& \ddots & & & &\ddots\\
& & \frac{\widetilde{I}}{N} & & & & \frac{\eta_N(\wo)}{N} A\\
\frac{\eta_1(\wo)}{N} A & & & -\frac{\widetilde{I}}{\nu N^2} &\cdots &\cdots & -\frac{\widetilde{I}}{\nu N^2}\\
& \ddots & &\vdots & \vdots & \vdots & \vdots\\
& & \frac{\eta_N(\wo)}{N} A\ & -\frac{\widetilde{I}}{\nu N^2} &\cdots &\cdots & -\frac{\widetilde{I}}{\nu N^2}
\end{pmatrix}
\begin{pmatrix}
\yb_1\\ \vdots \\\yb_N\\
\pb_1\\ \vdots \\ \pb_N
\end{pmatrix}=
\begin{pmatrix} \frac{\yb_d}{N}\\ \vdots \\
 \frac{\yb_d}{N} \\ \frac{\mathbf{f}}{N}\\ \vdots \\ \frac{\mathbf{f}}{N}
\end{pmatrix}.
\end{equation}
Next, let $\mathbf{z}=(\mathbf{z}_1,\dots,\mathbf{z}_{N_h})^\top \in \mathbb{R}^{(2N N_h)\times 1}$, where $\mathbf{z}_j=((\mathbf{y}_1)_j,\dots,(\mathbf{y}_N)_j,(\mathbf{p}_1)_j,\dots, (\mathbf{p}_N)_j))^\top \in \mathbb{R}^{2N\times 1}$. Notice that $\mathbf{z}_j$ corresponds to the application of $R_i$ to $\mathbf{x}$ (see \eqref{eq:restriction_action}), except for $u_i$ which has been previously eliminated.  
By reordering the unknowns as in $\mathbf{z}$, \eqref{eq:system_reduced} can be written as $S\mathbf{z}=\widetilde{\mathbf{b}}$ for a suitable $\widetilde{\mathbf{b}}$ and
\[S=\begin{pmatrix}
\widetilde{B} & B & \\
B & \widetilde{B} & B &\\
& B & \widetilde{B} & B &\\
 & & \ddots &  \ddots & \ddots &\\
& & & B &  \widetilde{B} & B &\\
& & & & B &  \widetilde{B}\\
\end{pmatrix}= \widetilde{I}\otimes \widetilde{B} + H\otimes B,\]
\[\widetilde{B}:=\begin{pmatrix}
\frac{I}{N} & D\\
D & -\frac{\unob\unob^\top}{\nu N^2}
\end{pmatrix},\quad B:=\begin{pmatrix}
0 & -\frac{D}{2}\\
-\frac{D}{2} & 0
\end{pmatrix},\quad H=\begin{pmatrix}
0 & 1\\
1 & 0 & 1\\
& \ddots &\ddots &\ddots \\
& &  1 & 0 &1\\
& &   & 1 & 0
\end{pmatrix},\]
where $I\in \mathbb{R}^{N\times N}$ is the identity matrix, $D$ is a diagonal matrix with $d_j:=\frac{2\eta_j(\wo)}{h^2N}$ on the diagonal, and $\unob=(1,\dots,1)\in \mathbb{R}^{N\times 1}$. In particular, a direct calcultation verifies that the iteration matrix of \eqref{eq:collective_smoothing} with $\theta=1$ and with this new order of unknowns is equal to
\[\mathcal{G}= \mathcal{I} - (\widetilde{I}\otimes \widetilde{B}^{-1})(\widetilde{I}\otimes \widetilde{B} + H\otimes B)=- H\otimes C,\]
with $C:=\widetilde{B}^{-1}B$, and $\mathcal{I}\in \mathbb{R}^{(2N_h N)\times (2N_h N)}$ being the identity matrix. We will next characterize precisely the spectrum of $\mathcal{G}$, which in turns gives an exact description of the convergence on the one-level collective smoother. To do so, we first study the spectrum of $C$ denoted by $\sigma(C)$.
\begin{lemma}[Spectrum of C]\label{lemma:spectrum_G}
The matrix $C$ has the spectrum
\[\sigma(C)=-\frac{1}{2}\left\{1,1-r\pm i\sqrt{(1-r)r}\right\},\]
with $r=\frac{\widehat{\mathbb{E}}\LQ \widetilde{\mathbf{d}}^{-2}\RQ}{\nu +\widehat{\mathbb{E}}\LQ \widetilde{\mathbf{d}}^{-2}\RQ}$, $\widetilde{\mathbf{d}}\in \mathbb{R}^{N\times 1}$, $(\widetilde{\mathbf{d}})_j=N d_j$, and $\widehat{\mathbb{E}}\LQ \widetilde{\mathbf{d}}^{-2}\RQ:=\frac{1}{N}\sum_{j=1}^N (\widetilde{\mathbf{d}})_j^{-2}$. The eigenvalue $\lambda=-\frac{1}{2}$ has algebraic multiplicity $2N-2$ and geometric multiplicity $N-1$.
\end{lemma}
\begin{proof}
Since $\frac{I}{N}$ and $D$ are non singular, to compute $C$ we use the exact formula for the inverse of $\widetilde{B}$. Setting $\Gamma:=\frac{1}{\nu N^2 +\unob^\top \frac{D^{-2}}{N}\unob}=\frac{1}{\nu N^2 +N^2 \widehat{\mathbb{E}}\LQ \widetilde{\mathbf{d}}^{-2}\RQ}$, with $(\widetilde{\mathbf{d}})_j=\frac{2\eta_j(\wo)}{h^2}$, we obtain
\begin{equation}\label{eq:C}
\begin{aligned}
C&=\widetilde{B}^{-1}B =\frac{1}{2}\begin{pmatrix}
-I +\frac{\Gamma}{N} D^{-1} \unob\unob^\top D^{-1} & -\Gamma D^{-1}\unob\unob^\top\\
\frac{D^{-1}}{N}-\frac{\Gamma}{N^2}D^{-2}\unob\unob^\top D^{-1} & -I +\frac{\Gamma}{N}D^{-2}\unob\unob^\top \end{pmatrix}\\
&=\frac{1}{2}\begin{pmatrix}
-I & 0\\
\frac{D^{-1}}{N} & -I
\end{pmatrix}+\frac{\Gamma N}{2}\begin{pmatrix}
\widetilde{\mathbf{d}}^{-1}\widetilde{\mathbf{d}}^{-\top} & -\widetilde{\mathbf{d}}^{-1}\unob^\top\\
-\widetilde{\mathbf{d}}^{-2}\widetilde{\mathbf{d}}^{-\top} & \widetilde{\mathbf{d}}^{-2}\unob^\top
\end{pmatrix}.
\end{aligned}
\end{equation}
For simplicity, we focus on $\widehat{C}:=-2C$, which can be written as
\[\widehat{C}=\underbrace{\begin{pmatrix}
I & 0\\
-\frac{D^{-1}}{N} & I
\end{pmatrix}}_{L}+ \mathbf{a}\mathbf{c}^\top,\quad \text{with}\quad \mathbf{a}:=\Gamma N  \begin{pmatrix}
-\widetilde{\mathbf{d}}^{-1}\\ \widetilde{\mathbf{d}}^{-2}
\end{pmatrix},\; \mathbf{c}:=\begin{pmatrix}
\widetilde{\mathbf{d}}^{-1}\\
-\unob
\end{pmatrix},\]
that is, $\widehat{C}$ is the sum of a lower triangular matrix plus a rank-one perturbation. Notice that $L$ has eigenvalue $\lambda=1$ with algebraic multiplicity $2N$ and geometric multiplicity $N$. The eigenspace associated to $\lambda=1$ is $E_{\lambda=1}(L):=\text{span}\left\{\mathbf{e}_j,\; j=N+1,\dots,2N\right\}$, $\mathbf{e}_j$ being the $j$-th canonical vector. Next, if $N> 2$, $\widehat{C}$ has still eigenvalue $\lambda=1$ since
for any vector $\mathbf{v}=(0,\mathbf{v}_2)$, $\mathbf{v}_2\in \mathbb{R}^{N\times 1}$, such that $\unob^\top \mathbf{v}_2=0$, we have
\[(L+\mathbf{a}\mathbf{c}^\top) \mathbf{v}= L\mathbf{v}=\mathbf{v}.\]
Therefore, $\lambda=1$ is an eigenvalue of $\widehat{C}$ with geometric multiplicity $N-1$.
To find the remaining eigenvalues, we take a $\lambda\neq 1$ and consider
\begin{equation}\label{eq:determinant}
\begin{aligned}
\det(L-\lambda I_{2N\times 2N} +\mathbf{a}\mathbf{c}^\top)&=\det(L-\lambda I_{2N\times 2N})\det(I_{2N\times 2N}+(L-\lambda I_{2N\times 2N})^{-1}\mathbf{a}\mathbf{c}^\top)\\
&=(1-\lambda)^{2N}\left(1+\mathbf{c}^\top (L-\lambda I_{2N\times 2N})^{-1} \mathbf{a}\right).
\end{aligned}
\end{equation}
A direct calculation leads to
\begin{equation}
\begin{aligned}
\mathbf{c}^\top (L-\lambda I_{2N\times 2N})^{-1} \mathbf{a}=(\mathbf{c}_1,\mathbf{c}_2)^\top \begin{pmatrix}
\frac{I}{1-\lambda} & 0 \\
\frac{D^{-1}}{N(1-\lambda)^2} &\frac{I}{1-\lambda}
\end{pmatrix}\begin{pmatrix}
\mathbf{a}_1\\ \mathbf{a}_2
\end{pmatrix}
\end{aligned},
\end{equation}
so that
\[\det(L-\lambda I_{2N\times 2N} +\mathbf{a}\mathbf{c}^\top)=(1-\lambda)^{2N-2}\left(\lambda^2-(2+\mathbf{a}^\top\mathbf{c})\lambda +1 +\mathbf{a}^\top\mathbf{c} +\mathbf{c}_2^\top \frac{D^{-1}}{N} \mathbf{a}_1\right),\]
from which we conclude that $\lambda=1$ has algebraic multiplicity $2(N-1)$.
The remaining eigenvalues must be solutions of the second order equation.
Using $\mathbf{a}^\top\mathbf{c}=-2\Gamma N \sum_{i=j}^N \widetilde{d}_j^{-2}$, 
$\mathbf{c}_2^\top \frac{D^{-1}}{N} \mathbf{a}_1=\Gamma N\sum_{j=1}^N \widetilde{d}_j^{-2}$, recalling the definition of $\Gamma$ and $r$, and dividing by $-\frac{1}{2}$, one obtains the solutions $\lambda_{2N-1,2N}=-\frac{1}{2}\left\{1-r\pm i \sqrt{(1-r)r}\right\}$, and the claim follows.
\end{proof}
\begin{remark}[Dependence on the regularization parameter]
The regularization parameter $\nu$ enters into our convergence analysis only in the definition of $r$. In particular as $\nu\rightarrow 0$, $r\rightarrow 1$ and $\lvert\lambda_{2N-1,2N}\rvert\rightarrow 0$, and the convergence of the collective multigrid does not deteriorate (see Lemma \ref{lemma:spectrum_G}). The robustness of the algorithm with respect to the (often troublesome) $\nu\rightarrow 0$ limit will be observed in the numerical experiments.
\end{remark}
From Lemma \ref{lemma:spectrum_G}, we deduce that $C$ admits the Jordan decomposition $CV=VJ$, with
\begin{equation}\label{eq:definition_J_V}
\begin{aligned}
J&=\begin{pmatrix}
-0.5 & 1\\
& -0.5 & \\
& & -0.5 &1 \\
& & &  -0.5\\
& & & &\ddots &\ddots\\
& & &  & &\lambda_{2N-1}\\
& & &  & & &&\lambda_{2N}
\end{pmatrix},\\
V&=[\mathbf{v}_1,\widehat{\mathbf{v}}_1,\mathbf{v}_{2},\widehat{\mathbf{v}}_2,\dots,\mathbf{v}_{2N-1},\mathbf{v}_{2N}],
\end{aligned}
\end{equation}
where $\mathbf{v}_j$, $j=1,\dots,N-1$, are the eigenvectors of $C$, $\widehat{\mathbf{v}}_j$, $j=1,\dots,N-1$, are the generalized eigenvectors satisfying $(C-\lambda_jI)\widehat{\mathbf{v}}_j=\mathbf{v}_j$, and $\mathbf{v}_{2N-1}$ and $\mathbf{v}_{2N}$ are the eigenvectors associated to the two remaining eigenvalues $\lambda_{2N-1,2N}$.

Exploiting the Kronecker structure of $\mathcal{G}$, we obtain immediately the following two corollaries.
\begin{corollary}[Similarity transformation of $\mathcal{G}$]\label{eq:Jordan_G}
For $i=1,\dots,N_h$ and $j=1,\dots 2N$, let $\delta_{i,j}:=-\mu_j\lambda_i$, where
$\lambda_i$ is an eigenvalue of $C$, and $\mu_j=2\cos\left(\frac{j\pi}{N_h+1}\right).$
Then, $\mathcal{G}$ satisfies $\mathcal{G} Y=Y\widetilde{J}$, where $\widetilde{J}$ is an upper triangular matrix with $\delta_{i,j}$ on the diagonal, and the $k$-th column of $Y$, with $k=i+j-1$ for some $i$ and $j$, is $Y_k=\bm{\varphi}_j\otimes V_i$, $V_i$ being the $i$-th column of $V$ defined in \eqref{eq:definition_J_V}, and $(\bm{\varphi}_j)_i:=\sin\left(\frac{ij\pi}{N_h+1}\right)$.
\end{corollary}
\begin{proof}
We first notice that $H$ is a tridiagonal Toeplitz matrix, and it is well-known (see \cite{https://doi.org/10.1002/nla.1811}) that has eigenvalues $\mu_j=2\cos\left(\frac{j\pi}{N_h+1}\right)$ and eigenvectors of the form $(\bm{\varphi}_j)_i=\sin\left(\frac{ij\pi}{N_h+1}\right)$.
Due to the properties of the Kronecker product, it is trivial to verify that
\[\mathcal{G}(\bm{\varphi}_j\otimes \mathbf{v}_i)=-(H\bm{\varphi}_j)\otimes (C\mathbf{v}_i)=-\mu_j\lambda_i (\bm{\varphi}_j\otimes \mathbf{v}_i).\]
If instead we consider a generalized eigenvector $\widehat{\mathbf{v}}_i$, using the Jordan decomposition, we have
\[\mathcal{G}(\bm{\varphi}_j\otimes \widehat{\mathbf{v}}_i)=-(H\bm{\varphi}_j)\otimes (C\widehat{\mathbf{v}}_i)=-\mu_j \lambda_i (\bm{\varphi}_j \otimes \widehat{\mathbf{v}}_i) -\mu_j (\bm{\varphi}_j \otimes \mathbf{v}_i),\]
and the claim follows.
\end{proof}
\begin{corollary}[Spectral radius of $\mathcal{G}$]\label{eq:spectrum_G}
The spectral radius of $\mathcal{G}$ is strictly smaller than $1$, and satisfies $\rho(\mathcal{G})\leq 1-\mathcal{O}\left(\frac{1}{N_h^2}\right)$. Therefore, the collective smoothing iteration converges.
\end{corollary}
\begin{proof}
Corollary \ref{eq:Jordan_G} shows that $\mathcal{G}$ is similar to the upper triangular matrix $\widetilde{J}$. Thus, its eigenvalues are equal to $\delta_{i,j}=-\mu_j\lambda_i$.
Observing that $\lvert \mu_j \rvert <2\lvert\cos\left(\frac{\pi}{N_h+1}\right)\rvert$ and $\lvert \lambda_i\rvert\leq 0.5$ for any $j, i$, the claim follows.
\end{proof}
\begin{remark}[Damping]
The analysis has been carried out for the relaxation parameter $\theta=1$. It is trivial to consider $\theta\neq 1$, since the iteration matrix is then
$\mathcal{G}_{\theta}:=(1-\theta) \mathcal{I} +\theta \mathcal{G}$.
\end{remark}
We next study the spectrum of the two-level collective multigrid algorithm, and assume that $N_{h}=2^{\ell}-1$ and $N^C_{h}=2^{\ell-1}-1$ for a $\ell\in\mathbb{N}$.
As maps between the fine and coarse meshes, we choose the full weighting restriction matrix,
\[\widetilde{R}:=\begin{pmatrix}
\frac{1}{2} & 1 &\frac{1}{2} \\
& & \frac{1}{2} & 1 & \frac{1}{2}\\
& & & \cdots \\
& & & & \frac{1}{2} & 1 &\frac{1}{2}\\\end{pmatrix}\in \mathbb{R}^{N_h^C\times N_h},
\]
and the linear interpolation operator $\widetilde{P}:=2 \widetilde{R}^\top$. In particular, the action of $\widetilde{R}$ and $\widetilde{P}$ on the frequencies $\bm{\varphi}_j$ can be characterized rigourously (see, e.g., \cite[Lemma 4.17]{ciaramella2022iterative}). Let $\bm{\phi}_j\in \mathbb{R}^{N_h^C\times 1}$ with $(\bm{\phi}_j)_i=\sin\left(\frac{2ij\pi}{N_h+1}\right)$, $j=1,\dots,N_h$ and $i=1,\dots,N_h^C$. Further define $c_j:=\cos\left(\frac{j\pi}{2(N_h+1)}\right)$ and $s_j:=\sin\left(\frac{j\pi}{2(N_h+1)}\right)$.
Then, for any $e_j,\; e_{\widetilde{j}}\in \mathbb{R}$, with $\widetilde{j}:=N_h+1-j$ and $j=1,\dots,\frac{N_h+1}{2}-1$, 
\begin{equation}\label{eq:action_restriction_prolongation}
\begin{aligned}
\widetilde{R} \begin{pmatrix}
\bm{\varphi}_j & \bm{\varphi}_{\widetilde{j}} 
\end{pmatrix}\begin{pmatrix}
e_j\\e_{\widetilde{j}}
\end{pmatrix}&=\widetilde{R}\left(e_j \bm{\varphi}_j+ e_{\widetilde{j}} \bm{\varphi}_{\widetilde{j}}\right)=(e_jc_j^2-e_{\widetilde{j}}s_j^2)\bm{\phi}_j=\bm{\phi}_j\begin{pmatrix}
c_j^2 & -s_j^2
\end{pmatrix}\begin{pmatrix}
e_j\\e_{\widetilde{j}}
\end{pmatrix},\\
\widetilde{P}\bm{\phi}_j&=(c_j^2\bm{\varphi}_j -s_j^2\bm{\varphi}_{\widetilde{j}})=\begin{pmatrix}
\bm{\varphi}_j & \bm{\varphi}_{\widetilde{j}}
\end{pmatrix}\begin{pmatrix}
c_j^2\\-s_j^2
\end{pmatrix}
\end{aligned}.
\end{equation}
Furthermore, $R\bm{\varphi}_{\jb}=0$ for $\jb:=\frac{N_h+1}{2}$.
The iteration matrix of the two-level algorithm with one-step of presmoothing and no post-smoothing is
\[T:=(I-RS_c^{-1}PS)\mathcal{G},\]
where $R=\widetilde{R}\otimes I$, $P=\widetilde{P}\otimes I$, and $S_C=RSP$.
\begin{lemma}\label{lemma:diagonalization_T}
The two-level operator $T$ is similar to a block diagonal matrix whose diagonal blocks are:
\begin{itemize}
\item[1] The matrices $T_{ji}:=\mathcal{G}_{ji}-R_{j}^\top\Pi_{ji}^{-1}R_{j}S_{ji}\mathcal{G}_{ji}\in \mathbb{R}^{4\times 4}$ for $j=1,\dots,\frac{N_h+1}{2}-1$ and $i=1,\dots,N-1$, with
\begin{equation}\label{eq:first_block}
\begin{aligned}
\mathcal{G}_{ji}&:= \begin{pmatrix}
\delta_{ji} &  & -\mu_j\\
& \delta_{\jt i} & &-\mu_{\jt}\\
& &  \delta_{ji} \\
& & & \delta_{\jt i}
\end{pmatrix},\quad S_{ji}:=\begin{pmatrix}
(1-\delta_{ji}) &  & -\mu_j\\
& (1-\delta_{\jt i}) & &-\mu_{\jt}\\
& &  (1-\delta_{ji}) \\
& & & (1-\delta_{\jt i})
\end{pmatrix},\\
R_j&:=\begin{pmatrix}
c_j^2& -s_j^2 \\
& & c_j^2 & -s_j^2
\end{pmatrix},\quad P_j=R_j^\top,\quad \Pi_{ji}:=R_j S_{ji}R_j^\top. 
\end{aligned}
\end{equation}
\item[2] The matrices $\mathcal{G}_{\jb i}=\begin{pmatrix}
\delta_{\jb i} & -\mu_{\jb}\\
 & \delta_{\jb i}
\end{pmatrix}\in \mathbb{R}^{2\times 2}$ for $\jb=\frac{N_h+1}{2}$, and $i=1,\dots,N-1$.
\item[3] The matrices $\widehat{T}_{ji}:=\widehat{\mathcal{G}}_{ji}-\widehat{R}_{j}^\top\widehat{\Pi}_{ij}^{-1}\widehat{R}_{j}\widehat{S}_{ji}\widehat{\mathcal{G}}_{ij}\in \mathbb{R}^{2\times 2}$ for $j=1,\dots,\frac{N_h+1}{2}-1$ and $i=2N-1, 2N$, with
\begin{equation}\label{eq:second_block}
\begin{aligned}
\widehat{\mathcal{G}}_{ji}&:= \begin{pmatrix}
\delta_{ji} & \\
& \delta_{\jt i}
\end{pmatrix},\quad \widehat{S}_{ji}:=\begin{pmatrix}
(1-\delta_{ji}) &\\
& (1-\delta_{\jt i})
\end{pmatrix},\\
\widehat{R}_j&:=\begin{pmatrix}
c_j^2& -s_j^2 \\
\end{pmatrix},\quad \widehat{P}_j=\widehat{R}_j^\top,\quad \widehat{\Pi}_{ji}:=c_j^4(1-\delta_{ji})+s_j^4(1-\delta_{ji}). 
\end{aligned}
\end{equation}
\item[4] The matrices $\widehat{\mathcal{G}}_{\jb i}=\begin{pmatrix}
\delta_{\jb i} & \\
 & \delta_{\jb i}
\end{pmatrix}\in \mathbb{R}^{2\times 2}$ for $\jb=\frac{N_h+1}{2}$, and $i=2N-1, 2N$.
\end{itemize}
\end{lemma}
\begin{proof}
The proof follows closely the arguments presented in \cite{ciaramella2022substructured,ciaramella2022spectral} for the study of two-level iterative methods. It consists in studying the action of $T$ onto suitably defined subspaces, showing that these subspaces are invariant, and finally deriving a matrix representation of $T$ into a new basis.
We start with the four dimensional subspaces $\mathcal{V}_{ji}:=\text{span}\left\{\bm{\varphi}_j\otimes \vb_i,\;
\bm{\varphi}_{\jt} \otimes \vb_i,\; \bm{\varphi}_{j} \otimes \widehat{\vb}_i,\; \bm{\varphi}_{\jt} \otimes \widehat{\vb}_i\right\}$, for $j=1,\dots,\frac{N_h+1}{2}-1$, $i=1,\dots,N-1$. For any quadruple of real numbers $e_j, e_{\jt},\widehat{e}_{j},\widehat{e}_{\jt}$, using $H\bm{\varphi}_j=\mu_j\bm{\varphi}_j$ and the Jordan decomposition of $C$, we obtain
\begin{equation}\label{eq:action_G}
\begin{aligned}&\mathcal{G}\begin{pmatrix}
\bm{\varphi}_j\otimes \vb_i & \bm{\varphi}_{\jt} \otimes \vb_i & \bm{\varphi}_j\otimes \widehat{\vb}_i & \bm{\varphi}_{\jt}\otimes \widehat{\vb}_i
\end{pmatrix}\begin{pmatrix}
e_j\\ e_{\jt}\\ \widehat{e}_{j} \\ \widehat{e}_{\jt}
\end{pmatrix}\\
&=\begin{pmatrix}
\bm{\varphi}_j\otimes \vb_i & \bm{\varphi}_{\jt} \otimes \vb_i & \bm{\varphi}_j\otimes \widehat{\vb}_i & \bm{\varphi}_{\jt}\otimes \widehat{\vb}_i
\end{pmatrix}
\begin{pmatrix}
\delta_{ji} &  & -\mu_j\\
& \delta_{\jt i} & &-\mu_{\jt}\\
& &  \delta_{ji} \\
& & & \delta_{\jt i}
\end{pmatrix}
\begin{pmatrix}
e_j\\ e_{\jt}\\ \widehat{e}_{j} \\ \widehat{e}_{\jt}
\end{pmatrix}.
\end{aligned}
\end{equation}
Next, since $\widehat{\vb}_i$ satisfies $(C-\lambda_i I)\widehat{\vb}_i=\vb_i$, it holds $B\widehat{\vb}_i=\widetilde{B}(\vb_i+\lambda_i\widehat{\vb}_i)$, hence,
\begin{equation}\label{eq:action_S}
\begin{aligned}&S\mathcal{G}\begin{pmatrix}
\bm{\varphi}_j\otimes \vb_i & \bm{\varphi}_{\jt} \otimes \vb_i & \bm{\varphi}_j\otimes \widehat{\vb}_i & \bm{\varphi}_{\jt}\otimes \widehat{\vb}_i
\end{pmatrix}\begin{pmatrix}
e_j\\ e_{\jt}\\ \widehat{e}_{j} \\ \widehat{e}_{\jt}
\end{pmatrix}\\
&=\begin{pmatrix}
\bm{\varphi}_j\otimes \vb_i & \bm{\varphi}_{\jt} \otimes \vb_i & \bm{\varphi}_j\otimes \widehat{\vb}_i & \bm{\varphi}_{\jt}\otimes \widehat{\vb}_i
\end{pmatrix}
\begin{pmatrix}
(1-\delta_{ji}) &  & -\mu_j\\
& (1-\delta_{\jt i}) & &-\mu_{\jt}\\
& &  (1-\delta_{ji}) \\
& & & (1-\delta_{\jt i})
\end{pmatrix}G_{ji}
\begin{pmatrix}
e_j\\ e_{\jt}\\ \widehat{e}_{j} \\ \widehat{e}_{\jt}
\end{pmatrix},
\end{aligned}
\end{equation}
and recalling \eqref{eq:action_restriction_prolongation},
\begin{equation}\label{eq:action_RSG}
\begin{aligned}
&RS\mathcal{G}\begin{pmatrix}
\bm{\varphi}_j\otimes \vb_i & \bm{\varphi}_{\jt} \otimes \vb_i & \bm{\varphi}_j\otimes \widehat{\vb}_i & \bm{\varphi}_{\jt}\otimes \widehat{\vb}_i
\end{pmatrix}\begin{pmatrix}
e_j\\ e_{\jt}\\ \widehat{e}_{j} \\ \widehat{e}_{\jt}
\end{pmatrix}\\
&=\begin{pmatrix}
\bm{\phi}_j\otimes \vb_i & \bm{\phi}_j\otimes \widehat{\vb}_i
\end{pmatrix}
\begin{pmatrix}
c_j^2& -s_j^2 \\
& & c_j^2 & -s_j^2
\end{pmatrix}S_{ji}G_{ji}
\begin{pmatrix}
e_j\\ e_{\jt}\\ \widehat{e}_{j} \\ \widehat{e}_{\jt}
\end{pmatrix}.
\end{aligned}
\end{equation}
We now consider the coarse correction.
\begin{equation}\label{eq:application_coarse}
\begin{aligned}
&S_c \begin{pmatrix} \bm{\phi}_j\otimes \vb_i &  \bm{\phi}_j\otimes \widehat{\vb}_i\end{pmatrix}\begin{pmatrix}
e^c_j\\ \widehat{e}^c_j
\end{pmatrix} =RSP \begin{pmatrix} \bm{\phi}_j\otimes \vb_i &  \bm{\phi}_j\otimes \widehat{\vb}_i\end{pmatrix}\begin{pmatrix}
e^c_j\\ \widehat{e}^c_j
\end{pmatrix}\\
&=RS \begin{pmatrix}
\bm{\varphi}_j\otimes \vb_i & \bm{\varphi}_{\jt} \otimes \vb_i & \bm{\varphi}_j\otimes \widehat{\vb}_i & \bm{\varphi}_{\jt}\otimes \widehat{\vb}_i
\end{pmatrix}R_{j}^\top \begin{pmatrix}
e^c_j\\ \widehat{e}^c_j
\end{pmatrix}\\
&=\begin{pmatrix} \bm{\phi}_j\otimes \vb_i &  \bm{\phi}_j\otimes \widehat{\vb}_i\end{pmatrix}R_{j}S_{ji} R_{j}^\top \begin{pmatrix}
e^c_j\\ \widehat{e}^c_j,
\end{pmatrix}=\begin{pmatrix} \bm{\phi}_j\otimes \vb_i &  \bm{\phi}_j\otimes \widehat{\vb}_i\end{pmatrix}  \Pi_{ji} \begin{pmatrix}
e^c_j\\ \widehat{e}^c_j,
\end{pmatrix}
\end{aligned}
\end{equation}
which implies
\[S_c^{-1}\begin{pmatrix} \bm{\phi}_j\otimes \vb_i &  \bm{\phi}_j\otimes \widehat{\vb}_i\end{pmatrix}=\begin{pmatrix} \bm{\phi}_j\otimes \vb_i &  \bm{\phi}_j\otimes \widehat{\vb}_i\end{pmatrix}\Pi_{ij}^{-1}.\]
Putting all together, we get
\begin{equation}\label{eq:actionT}
\begin{aligned}
&T\begin{pmatrix}
\bm{\varphi}_j\otimes \vb_i & \bm{\varphi}_{\jt} \otimes \vb_i & \bm{\varphi}_j\otimes \widehat{\vb}_i & \bm{\varphi}_{\jt}\otimes \widehat{\vb}_i
\end{pmatrix}\begin{pmatrix}
e_j\\ e_{\jt}\\ \widehat{e}_{j} \\ \widehat{e}_{\jt}
\end{pmatrix}\\
&=\begin{pmatrix}
\bm{\varphi}_j\otimes \vb_i & \bm{\varphi}_{\jt} \otimes \vb_i & \bm{\varphi}_j\otimes \widehat{\vb}_i & \bm{\varphi}_{\jt}\otimes \widehat{\vb}_i
\end{pmatrix}\underbrace{\left(\mathcal{G}_{ji}-R_{j}^\top\Pi_{ji}^{-1}R_{j}S_{ji}\mathcal{G}_{ji}\right)}_{T_{ji}} \begin{pmatrix}
e_j\\ e_{\jt}\\ \widehat{e}_{j} \\ \widehat{e}_{\jt}
\end{pmatrix}.
\end{aligned}
\end{equation}
This conclude the first part of the proof. We now consider the subspaces spanned by
$\bm{\varphi}_{\jb} \otimes \vb_i$, $\bm{\varphi}_{\jb} \otimes \widehat{\vb}_i$ for $i=1,\dots,N-1$.
Since $R\bm{\varphi}_{\jb}=0$, we immediately have
\[T \begin{pmatrix}
\bm{\varphi}_{\jb}\otimes \vb_i & \bm{\varphi}_{\jb}\otimes \widehat{\vb}_i \end{pmatrix}
\begin{pmatrix}
e_{\jb}\\ \widehat{e}_{\jb} 
\end{pmatrix}=\begin{pmatrix}
\bm{\varphi}_{\jb}\otimes \vb_i & \bm{\varphi}_{\jb}\otimes \widehat{\vb}_i \end{pmatrix} \mathcal{G}_{\jb i}
\begin{pmatrix}
e_{\jb}\\ \widehat{e}_{\jb} 
\end{pmatrix},\quad  \mathcal{G}_{\jb i}:=\begin{pmatrix}
\delta_{\jb i} &-\mu_{\jb}\\
& \delta_{\jb i}
\end{pmatrix},\]
and this proves the second claim.
As third set of subspaces, we consider those spanned by respectively $(\bm{\varphi}_j\otimes v_{2N-1},\bm{\varphi}_{\jt}\otimes v_{2N-1})$, and $(\bm{\varphi}_j\otimes v_{2N},\bm{\varphi}_{\jt}\otimes v_{2N})$.
Following the same calculations of the first part of the proof we obtain for $i=2N-1$ and $i=2N$,
\begin{equation}
\begin{aligned}
&T\begin{pmatrix}
\bm{\varphi}_j\otimes v_{i} & \bm{\varphi}_{\jt}\otimes v_{i}
\end{pmatrix}\begin{pmatrix}
e_{j} \\ e_{\jb}
\end{pmatrix}\\
&=\begin{pmatrix}
\bm{\varphi}_j\otimes v_{i} & \bm{\varphi}_{\jt}\otimes v_{i}
\end{pmatrix}
\left(\widehat{\mathcal{G}}_{ji}-\widehat{R}^\top_j \widehat{\Pi}_{ji}^{-1} \widehat{R}_j\widehat{S}_{ji}\mathcal{G}_{ji}\right)\begin{pmatrix}
e_{j} \\ e_{\jb}
\end{pmatrix}
\end{aligned}
\end{equation}
The proof of the fourth claim is identical to that of the second part and it is skipped for the sake of brevity.
By considering a matrix $V$ that has column-block wise the basis for the subspaces we considered, it is immediate to deduce that $TV=V\widetilde{T}$, where $\widetilde{T}$ is a block diagonal matrix with the blocks we computed.
\end{proof}
\begin{remark}[Generalization to arbitrary pre- and post-smoothing steps]
Lemma \eqref{lemma:diagonalization_T} can be readily generalized to cover $n_1$ pre-smoothing steps and $n_2$ post-smoothing steps, but taking suitable powers of the matrices $\mathcal{G}_{ji}$, $\widehat{\mathcal{G}}_{\jb i}$, $\widehat{\mathcal{G}}_{ji}$ and $\widehat{\mathcal{G}}_{\jb i}$. For instance, the matrix $T_{ji}$ of part one becomes
\[T_{ji}:=\mathcal{G}_{ji}^{n_2}(I_{4\times 4}-R_j^\top\Pi_{ji}^{-1}R_jS_{ji})\mathcal{G}_{ji}^{n_1}.\]
\end{remark}

\begin{figure}
\centering
\includegraphics[scale=0.38]{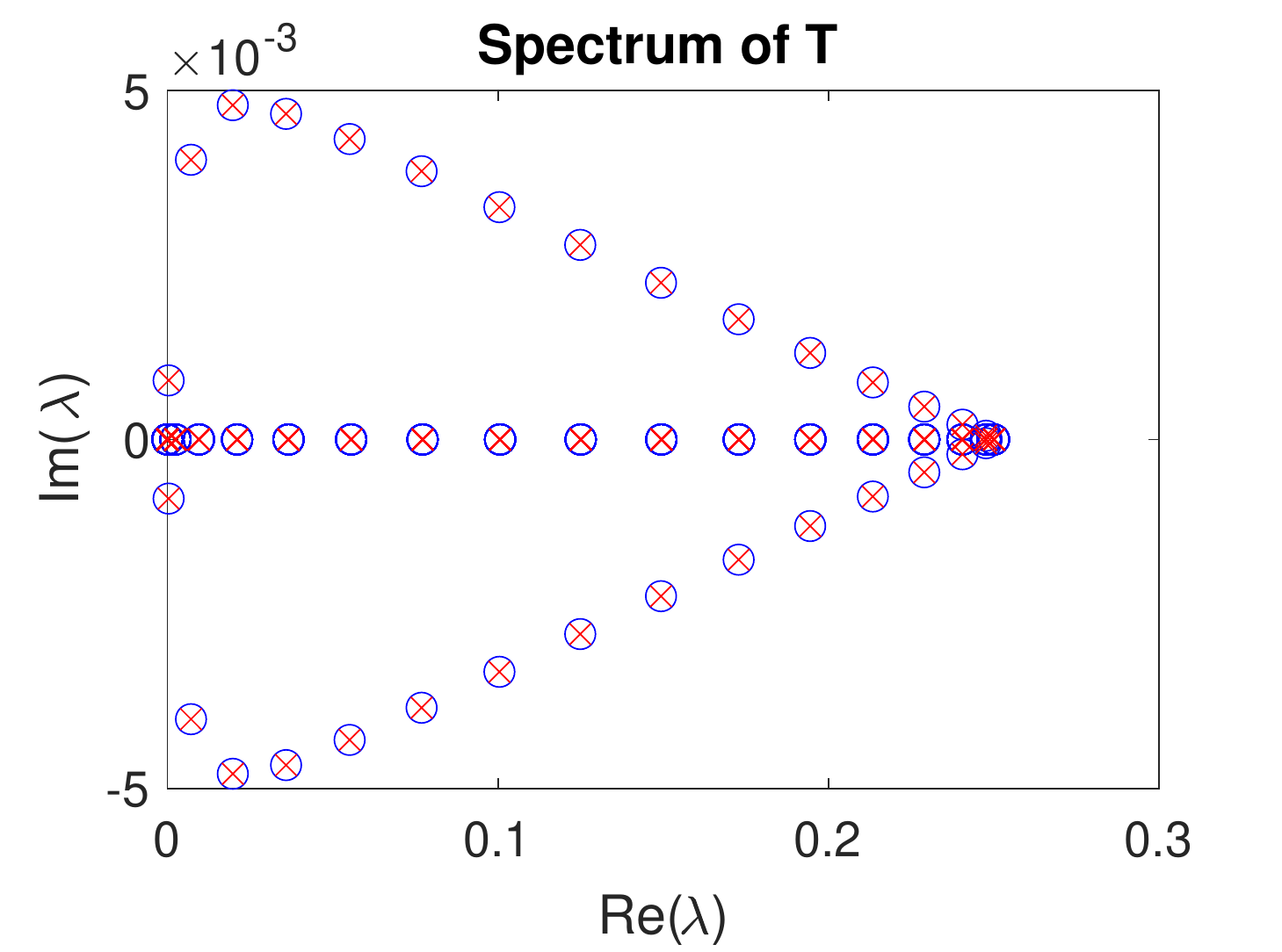}
\includegraphics[scale=0.38]{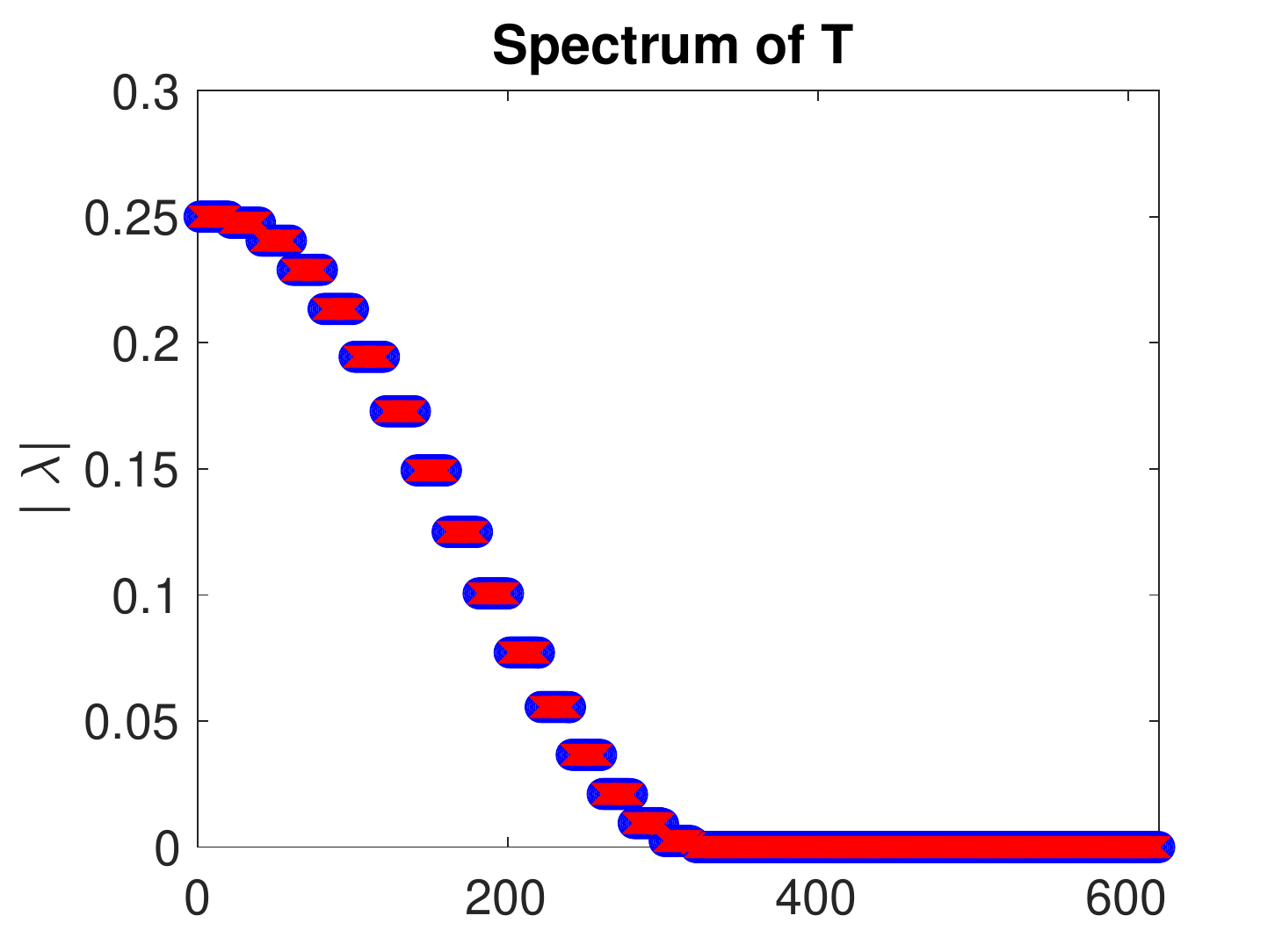}
\includegraphics[scale=0.38]{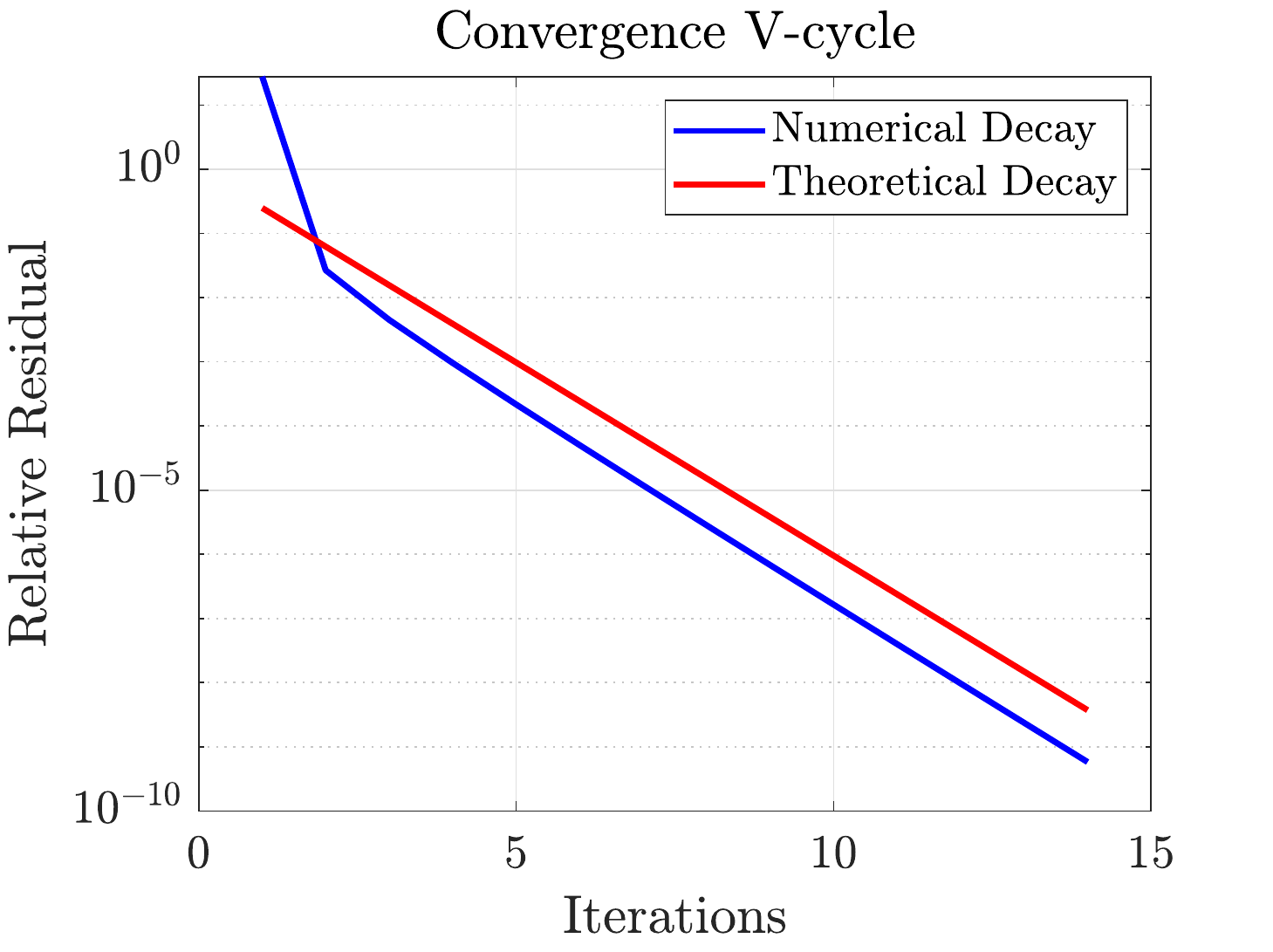}
\caption{Top row: Graphical representation of the spectrum of $T$ for $N_h=31$, $N=10$, $\nu=10^{-2}$ and $n_1=n_2=1$. The blue circles are obtained by computing numerically the eigenvalues of T. The red crosses are obtained through the formulae of Theorem \ref{theorem:spectrum_two_level}. Bottom row: comparison between the numerical and theoretical convergence of the two-level algorithm.}\label{fig:spectrum_T}
\end{figure}

\begin{theorem}[Spectrum and convergence of the two-level algorithm]\label{theorem:spectrum_two_level}
The spectrum of the matrix $T=\mathcal{G}^{n_2}(I-R S_c^{-1}PS)\mathcal{G}^{n_1}$ is
\begin{equation}\label{eq:spectrum_T}
\begin{aligned}
\sigma(T)&=\left\{0\right\}\cup\\
&\left\{\frac{c_j^4(1-\delta_{ji})\delta_{\jt i}^{n_1+n_2}+s_j^4(1-\delta_{\jt i})\delta_{j i}^{n_1+n_2}}{c_j^4(1-\delta_{ji})+s_j^4(1-\delta_{\jt i})},\;j=1,\dots,\frac{N_h+1}{2}-1,\;i=1,\dots,2N\right\}.
\end{aligned}
\end{equation}
Further, the spectral radius of $T$ is strictly smaller than 1, hence the two-level collective multigrid algorithm converges.
\end{theorem}
\begin{proof}
Since $T$ is similar to a block diagonal matrix, with blocks defined in Lemma \ref{lemma:diagonalization_T}, it is sufficient to compute the spectrum of each block.
Further, the spectrum of $T$ is equal to that of $(I-R S_c^{-1}PS)\mathcal{G}^{n_1+n_2}$. Hence, we start considering the blocks $T_{ji}=(I_{4\times 4}-R_j^\top\Pi_{ji}^{-1}R_jS_{ji})\mathcal{G}_{ji}^{n_1+n_2}$. Direct calculations show that
\begin{equation}
\begin{aligned}
(I_{4\times 4}-R_j^\top\Pi_{ji}^{-1}R_jS_{ji})=
 \left(\begin{array}{@{}c|c@{}}
  \begin{matrix}
1- \frac{c_j^4(1-\delta_{ji})}{\gamma} & \frac{c_j^2s_j^2(1-\delta_{\jt i})}{\gamma} \\
  \frac{c_j^2s_j^2(1-\delta_{j i})}{\gamma} & 1-\frac{s_j^4(1-\delta_{\jt i})}{\gamma}
  \end{matrix}
  & \Large{X} \\
\hline
   &
  \begin{matrix}
1- \frac{c_j^4(1-\delta_{ji})}{\gamma} & \frac{c_j^2s_j^2(1-\delta_{\jt i})}{\gamma} \\
  \frac{c_j^2s_j^2(1-\delta_{j i})}{\gamma} & 1-\frac{s_j^4(1-\delta_{\jt i})}{\gamma}
  \end{matrix}
\end{array}\right)
\end{aligned},
\end{equation}
 where the expression of $X\in \mathbb{R}^{4\times 4}$ will not be needed in the following and $\gamma:=c_j^4(1-\delta_{ji})+s_j^4(1-\delta_{\jt i})$. Since the product of two upper triangular matrices is still upper triangular, it follows that
\[(I_{4\times 4}-R_j^\top\Pi_{ji}^{-1}R_jS_{ji})\mathcal{G}_{ji}^{n_1+n_2}=\begin{pmatrix}
K & \widetilde{X} \\
& K
\end{pmatrix},\]
with
\[K:=\frac{1}{\gamma}\begin{pmatrix}
s_j^4(1-\delta_{\jt i})\delta_{ji}^{n_1+n_2} & c_j^2s_j^2(1-\delta_{\jt i})\delta_{\jt i}^{n_1+n_2} \\
  c_j^2s_j^2(1-\delta_{j i})\delta_{j i}^{n_1+n_2} & c_j^4(1-\delta_{ji })\delta_{\jt i}^{n_1+n_2}
  \end{pmatrix},\]
and whose eigenvalues are $\kappa^{ji}_1=\frac{c_j^4(1-\delta_{ji})\delta_{\jt i}^{n_1+n_2}+s_j^4(1-\delta_{\jt i})\delta_{j i}^{n_1+n_2}}{c_j^4(1-\delta_{ji})+s_j^4(1-\delta_{\jt i})}$ and $\kappa_2=0$. Next, $\mathcal{G}^{n_1+n_2}_{\jb i}$ and $\widehat{\mathcal{G}}^{n_1+n_2}_{\jb i}$ have trivially eigenvalues equal to $\delta^{n_1+n_2}_{\jb i}$, which are all equal to zero since $\mu_{\bar{j}}=0$. Further, direct calculations show that $\widehat{T}_{ji}$ has also two eigenvalues equal, again, to $\kappa^{ji}_1$ and $\kappa_2$. Taking into account the range of the indices of $j$ and $i$ for each blocks, we obtain the characterization of the spectrum, and since $\lvert\delta_{ji}\rvert<1$ and $\lvert\delta_{\jt i}\rvert<1$, we conclude that the spectral radius of $T$ is smaller than one.
\end{proof}
Fig. \ref{fig:spectrum_T} shows the spectrum of $T$ where, for visualization purposes, we set $N_h=31$ and $N=10$. In particular, the right panel shows that the spectrum is grouped into $\frac{N_h-1}{2}$ clusters, in which each eigenvalue is repeated approximately $2N$ times (approximately, because $C$ has two eigenvalues, $\lambda_{2N-1,2N}$ slightly different from $0.5$.)
\begin{remark}[Extension of the analysis to the deterministic setting]
Our analysis also represents a novel approach to study the convergence of collective smoothing iterations in the case of a deterministic PDE constraint by setting $N=1$. Retracing the analysis, we observe that $C$ has only two eigenvalues equal to $\lambda_{2N-1,2N}$ and $\mathcal{G}$ is diagonal. $T$ can then be diagonalized more easily, and its spectrum is still characterized by \eqref{eq:spectrum_T}, where the index $i$ assumes only the values $2N-1$ and $2N$.
\end{remark}
\begin{remark}[Extension to the two and three dimensional physical space]
The analysis could be extended to square or cube domains. Due to the Kronecker product structure between spatial and probability quantities, only the matrix $H$ would have to change, and its eigenvectors would be the tensorized product of sine functions. Similarly, the action of the operators $\widetilde{R}$ and $\widetilde{P}$ would be represented by more complicated matrices.
\end{remark}

This concludes our theoretical study of the convergence of the two-level collective multigrid algorithm. The next sections will focus on analyzing its numerical performances in different cases.

\subsection{Numerical experiments}
We now show the performance of Alg. \ref{Alg:algorithm-Vcycle} and its robustness with respect to several parameters for the solution of \eqref{eq:full_space_optimality_system}. We first consider the state equation
\begin{equation}\label{eq:state_equation2}
a_{\omega}(y_{\omega},v)=\int_\D \kappa(x,\omega) \nabla y(x,\omega)\cdot\nabla v(x)\ dx=\int_\D u(x)v(x)\ dx,\quad \forall v\in V,\ \PP\text{-a.e. }\omega \in \Omega,
\end{equation} 
in the L-shaped domain $\D=(0,1)^2\setminus \overline{(0.5,1)}^2$ discretized with a regular mesh of squares of edge $h_{\ell}=2^{-\ell}$, which are then decomposed into two right triangles. We choose $\kappa(x,\omega)$ as an approximated log-normal diffusion field 
\begin{equation}\label{eq:log_normal}
\kappa(x,\omega)=e^{\sigma \sum_{j=1}^M \sqrt{\lambda_j}b_j(x)N_j(\omega)}\approx e^{g(x,\omega)},
\end{equation}
where $g(x,\omega)$ is a mean zero Gaussian field with Covariance function $Cov_g(x,y)=\sigma^2 e^{\frac{-\|x-y\|_{2}^2}{L^2}}$. The parameter $\sigma^2$ tunes the variance of the random field, while $L$ denotes the correlation length. The pairs $(b_j(x),\sigma^2\lambda_j)$ are the eigenpairs of $T:L^2(\D)\rightarrow L^2(\D)$, $(Tf)(x)=\int_\D Cov_g(x,y)f(y)\ dy$, and $N_j\stackrel{iid}{\sim} \mathcal{N}(0,1)$. Assumption \ref{ass:diff} is satisfied since $a_{\min}(\omega)=\left(\text{ess}\inf_{x\in \D} \kappa(x,\omega)\right)^{-1}$ and $a_{\max}(\omega)=\|\kappa(\cdot,\omega)\|_{L^\infty(\D)}$ are in $L^p(\Omega)$ for every $p<\infty$ \cite{Charrier}. The target state is $y_d=e^{y^2}\sin(2\pi x)\sin(2\pi y)$.

Table \ref{Tab:expected} shows the number of V-cycle iterations (Alg. \ref{Alg:algorithm-Vcycle}) and of GMRES iterations preconditioned by the V-cycle to solve \eqref{eq:full_space_optimality_system} up to a tolerance of $10^{-9}$ on the relative (unpreconditioned) residual.
Inside the V-cycle algorithm, we use $n_1=n_2=2$ pre- and post-smoothing iterations based on the Jacobi relaxation \eqref{eq:collective_smoothing} with a damping parameter $\theta=0.5$ (the same value will be used for all numerical experiments in this manuscript). Numerically, we observed that Gauss-Seidel relaxations lead to very similar results. The number of levels of the V-cycle hierachy is denoted with $N_L$. The size of the largest linear system solved per sub-table is denoted by $N_{\max}=(2N+1)N_h$.
\begin{table}\caption{Number of V-cycle (left) and preconditioned GMRES (right) iterations to solve \eqref{eq:full_space_optimality_system} for a linear quadratic problem on the L-shaped domain $\D=(0,1)^2\setminus \overline{(0.5,1)}^2$ with a distributed control.}\label{Tab:expected}
{\small
\centering
\begin{tabular}{| c | c | c | c | c|}
\hline
$\nu$ &  $ 10^{-2}$ & $10^{-4}$ &  $10^{-6}$ & $10^{-8}$\\ \hline\hline
It. & 18\; $\vert$ \;11 &  19\; $\vert$ \;13  &  19\; $\vert$ \;15 & 19\; $\vert$ \;15 \\
\hline
\end{tabular}\\
\centering
$N_h=705$, $N=125$, $N_L=3$, $\sigma^2=0.5$, $L^2=0.5$, $N_{\max}=1.77\;10^5$.\\\vspace*{0.2cm}
\begin{tabular}{| c | c | c | c | c|}
\hline
$\sigma^2$ &  0.1 & 0.5 &  1 & 1.5\\ \hline\hline
It. & 19 \; $\vert$ \; 13 & 19 \; $\vert$ \; 13  & 20 \; $\vert$ \; 13 & 20  \; $\vert$ \; 13\\
\hline
\end{tabular}\\
\centering
$N_h=705$, $N=125$, $N_L=3$, $\nu=10^{-4}$, $L^2=0.5$, $N_{\max}=1.77\;10^5$.\\\vspace*{0.2cm}
\centering
\begin{tabular}{| c | c | c | c |}
\hline
$N_h$($N_L$) & 161 (2) &  705 (3) & 2945 (4)\\ \hline\hline
It. & 19 \; $\vert$ \; 13 &  19 \; $\vert$ \; 13 & 20 \; $\vert$ \; 13 \\
\hline
\end{tabular}\\
\centering
$N=125$, $\nu=10^{-4}$, $\sigma^2=0.5$, $L^2=0.5$, $N_{\max}=7.39\;10^5$.\\\vspace{0.2cm}
\centering
\begin{tabular}{| c | c | c | c | c|}
\hline
$N$ & 8 & 27 & 64 & 125\\ \hline\hline
It. & 19 \; $\vert$ \; 13  &  19 \; $\vert$ \; 13 &  19 \; $\vert$ \; 13 & 19 \; $\vert$ \; 13 \\
\hline
\end{tabular}\\
\centering
$N_h=705$, $N_L=3$, $\nu=10^{-4}$, $\sigma^2=0.5$, $L^2=0.5$, $N_{\max}=1.77\;10^5.$\\
\centering
\begin{tabular}{| c | c | c | c | c|}
\hline
$N$ & 100 & 500 & 1000 &2000 \\ \hline\hline
It. & 22  \; $\vert$ \; 16 &  22  \; $\vert$ \; 15 & 22  \; $\vert$ \; 15 & 22  \; $\vert$ \; 15\\
\hline
\end{tabular}\\
\centering
$N_h=705$, $\nu=10^{-4}$, $N_L=3$, $\sigma^2=1.5$, $L^2=0.1$, $N_{\max}=2.82\;10^6$.\\\vspace{0.3cm}
}
\end{table}

\noindent The first four sub-tables are based on a discretization of the probability space using the Stochastic Collocation method \cite{babuvska2010stochastic} on Gauss-Hermite tensorized quadrature nodes, since for $L^2=0.5$, setting $M=3$ into \eqref{eq:log_normal} is enough to preserve $99\%$ of the variance. In the fifth sub-table we set $L^2=0.1$ and use the Monte Carlo method, since we need $M=15$ random variables to preserve $99\%$ of the variance of the random field, and the Stochastic Collocation method suffers the curse of dimensionality.
Remark that the multigrid algorithm is robust with respect to all parameters considered, namely the regularization parameter,  the variance of the random field, the number of levels as the fine grid is refined, and the number of samples to discretize the probability space.

We mention that a family of block diagonal preconditioners for saddle-point matrices such as \eqref{eq:matrix_saddle_point} were recently proposed in \cite{Kouri2018}. A detailed theoretical analysis was developed in \cite{nobile_vanzan} for distributed controls, in a more general setting than the one considered in this manuscript that covers a general finite element discretization of a $d$-dimensional domain, a general elliptic bilinear form, and an additional variance term in the cost functional. Their main attractive feature is the possibility to precondition fully in parallel the $2N$ PDEs. Nevertheless, their convergence deteriorates as $\nu\rightarrow 0$ (as several preconditioners built on the same technique see, e.g., \cite{rees2010optimal,pearson2012new}), so that these preconditioners are hardly effective when $\nu$ is smaller than, say, $10^{-3}/10^{-4}$. The robustness of the multigrid algorithm as $\nu\rightarrow 0$ is definitely one of its most interesting properties. In terms of mesh refinement, both approaches are robust, provided that the
$2N$ PDEs constraints are suitable preconditioned (e.g., with multigrid) in the approach of \cite{Kouri2018,nobile_vanzan}.
Concerning the refinement of the discretization of the probability space, both methods are robust, and interestingly, both convergence analyses show a dependence on the approximated expected value of the square inverse of the coercivity constants of the stiffness matrices. One current disadvantage of the multigrid algorithm is the lack of coarsening with respect to the number of samples $N$, since the solution of the coarse problem might represent a bottle neck for very fine discretizations. In these circumstances, the capability of \cite{Kouri2018,nobile_vanzan} to handle the PDE constraints in parallel may be beneficial.

Next, we consider the same problem \eqref{eq:state_equation2}-\eqref{eq:log_normal} posed in the unit square domain $\mathcal{D}=(0,1)^2$ with either a local control acting on the subset $\mathcal{D}_0=(0.25, 0.75)^2\subset \mathcal{D}$, or a Neumann boundary control acting on $\Gamma=(0,1)\times \left\{0\right\}\subset \partial \mathcal{D}$. Tables \ref{Tab:expected_local} and \ref{Tab:expected_boundary} report the performances of the multigrid algorithm for these two cases. We stress once more the excellent robustness and efficiency of the multigrid algorithm in all regimes.

\begin{table}\caption{Number of V-cycle (left) and preconditioned GMRES (right) iterations to solve \eqref{eq:full_space_optimality_system} for a linear quadratic problem on the square domain $\D=(0,1)^2$ with a local control acting on $\D_0=(0.25,0.75)^2$.}\label{Tab:expected_local}
{\small
\centering
\begin{tabular}{| c | c | c | c | c|}
\hline
$\nu$ &  $ 10^{-2}$ & $10^{-4}$ &  $10^{-6}$ & $10^{-8}$\\ \hline\hline
It. & 17\; $\vert$ \;11 &  20\; $\vert$ \;13  &  26\; $\vert$ \;16 & 26\; $\vert$ \;18 \\
\hline
\end{tabular}\\
\centering
$N_h=961$, $N=125$, $N_L=3$, $\sigma^2=0.5$, $L^2=0.5$.\\\vspace*{0.2cm}
\begin{tabular}{| c | c | c | c | c|}
\hline
$\sigma^2$ &  0.1 & 0.5 &  1 & 1.5\\ \hline\hline
It. & 20 \; $\vert$ \; 13 & 20 \; $\vert$ \; 13  & 20 \; $\vert$ \; 13 & 19  \; $\vert$ \; 13\\
\hline
\end{tabular}\\
\centering
$N_h=961$, $N=125$, $N_L=3$, $\nu=10^{-4}$, $L^2=0.5$.\\\vspace*{0.2cm}
\centering
\begin{tabular}{| c | c | c | c |}
\hline
$N_h$($N_L$) & 225 (2) &  961 (3) & 3969 (4)\\ \hline\hline
It. & 19 \; $\vert$ \; 12 &  20 \; $\vert$ \; 13 & 20 \; $\vert$ \; 13 \\
\hline
\end{tabular}\\
\centering
$N=125$, $\nu=10^{-4}$, $\sigma^2=0.5$, $L^2=0.5$.\\\vspace{0.2cm}
\centering
\begin{tabular}{| c | c | c | c | c|}
\hline
$N$ & 8 & 27 & 64 & 125\\ \hline\hline
It. & 20 \; $\vert$ \; 13  &  20 \; $\vert$ \; 13 &  20 \; $\vert$ \; 13 & 20 \; $\vert$ \; 13 \\
\hline
\end{tabular}\\
\centering
$N_h=961$, $N_L=3$, $\nu=10^{-4}$, $\sigma^2=0.5$, $L^2=0.5$.\\
\centering
\begin{tabular}{| c | c | c | c|}
\hline
$N$ & 100 & 1000 &2000 \\ \hline\hline
It. & 20  \; $\vert$ \; 14 &  21  \; $\vert$ \; 15 & 21  \; $\vert$ \; 14\\
\hline
\end{tabular}\\
\centering
$N_h=961$, $\nu=10^{-4}$, $N_L=3$, $\sigma^2=1.5$, $L^2=0.1$.\\\vspace{0.3cm}
}
\end{table}

\begin{table}\caption{Number of V-cycle (left) and preconditioned GMRES (right) iterations to solve \eqref{eq:full_space_optimality_system} for a linear quadratic problem on the square domain $\D=(0,1)^2$ with a boundary control acting on $\Gamma=(0,1)\times \left\{0\right\}$.}\label{Tab:expected_boundary}
{\small
\centering
\begin{tabular}{| c | c | c | c | c|}
\hline
$\nu$ &  $ 10^{-2}$ & $10^{-4}$ &  $10^{-6}$ & $10^{-8}$\\ \hline\hline
It. & 17\; $\vert$ \;14 &  22\; $\vert$ \;15  &  23\; $\vert$ \;17 & 21\; $\vert$ \;16 \\
\hline
\end{tabular}\\
\centering
$N_h=992$, $N=125$, $N_L=3$, $\sigma^2=0.5$, $L^2=0.5$.\\\vspace*{0.2cm}
\begin{tabular}{| c | c | c | c | c|}
\hline
$\sigma^2$ &  0.1 & 0.5 &  1 & 1.5\\ \hline\hline
It. & 16 \; $\vert$ \; 13 & 17 \; $\vert$ \; 14  & 17 \; $\vert$ \; 14 & 17  \; $\vert$ \; 14\\
\hline
\end{tabular}\\
\centering
$N_h=992$, $N=125$, $N_L=3$, $\nu=10^{-4}$, $L^2=0.5$.\\\vspace*{0.2cm}
\centering
\begin{tabular}{| c | c | c | c |}
\hline
$N_h$($N_L$) & 240 (2) &  992 (3) & 4032 (4)\\ \hline\hline
It. & 16 \; $\vert$ \; 12 &  17 \; $\vert$ \; 14 & 21 \; $\vert$ \; 16 \\
\hline
\end{tabular}\\
\centering
$N=125$, $\nu=10^{-4}$, $\sigma^2=0.5$, $L^2=0.5$.\\\vspace{0.2cm}
\centering
\begin{tabular}{| c | c | c | c | c|}
\hline
$N$ & 8 & 27 & 64 & 125\\ \hline\hline
It. & 16 \; $\vert$ \; 13  &  17 \; $\vert$ \; 14 &  17 \; $\vert$ \; 14 & 17 \; $\vert$ \; 14 \\
\hline
\end{tabular}\\
\centering
$N_h=992$, $N_L=3$, $\nu=10^{-4}$, $\sigma^2=0.5$, $L^2=0.5$.\\
\centering
\begin{tabular}{| c | c | c | c|}
\hline
$N$ & 100 & 1000 &2000 \\ \hline\hline
It. & 18  \; $\vert$ \; 15 &  19  \; $\vert$ \; 16 & 20  \; $\vert$ \; 16\\
\hline
\end{tabular}\\
\centering
$N_h=992$, $\nu=10^{-4}$, $N_L=3$, $\sigma^2=1.5$, $L^2=0.1$.\\\vspace{0.3cm}
}
\end{table}

\section{An optimal control problem under uncertainty with box-constraints and $L^1$ penalization}\label{Sec:l1}
In this section, we consider the nonsmooth OCPUU\footnote{To keep a light notation, we omitted the continuous embedding operator from $L^2(\Omega;V)$ to $L^2(\Omega;L^2(\D))$ and from $L^2(\D)$ to $V^\prime$, see Section \ref{Sec:quadratic} and, e.g., \cite[Section 2.13]{troltzsch2010optimal}.}
\begin{equation}\label{eq:l1_OCP}
\begin{aligned}
&\min_{u\in U_{ad}} \frac{1}{2}\E\LQ \|y_\omega(u)-y_d\|^2_{L^2(\D)}\RQ + \frac{\nu}{2}\|u\|^2_{L^2(\D)} + \beta \|u\|_{L^1(\D)},\\
&\quad \text{subject to}\\
&a_\omega(y_\omega(u),v)=(u+f,v)_{L^2(\D)},\quad \forall v \in V,\ \PP\text{-a-e. } \omega\in \Omega,\\
&U_{ad}:=\left\{v\in L^2(\D): a\leq u \leq b\quad \text{almost everywhere in }\D\right\},
\end{aligned}
\end{equation}
with $a<0<b$ and $\nu,\beta>0$.
Deterministic OCPs with a $L^1$ penalization lead to optimal controls which are sparse, i.e. they are nonzero only on certain regions of the domain $\D$ \cite{stadler2009elliptic,casas2017review}. Sparse controls can be of great interest in applications, because it is often not desirable, or even impossible, to control the system over the whole domain $\D$. For sparse OCPUU, we mention \cite{li2019sparse} where the authors considered both a simplified version of \eqref{eq:l1_OCP} in which the randomness enters linearly into the state equation as a force term, and a different optimization problem whose goal is to find a stochastic control $u(\omega)$ which has a similar sparsity pattern regardless of the realization $\omega$. Note further that the assumption $\nu>0$ does not eliminate the nonsmoothness of the objective functional, but it regularizes the optimal solution $u$, and is needed to use the fast optimation algorithm described in the following.

The well-posedness of \eqref{eq:l1_OCP} follows directly from standard variational arguments \cite{troltzsch2010optimal,hinze2008optimization}, being $U_{ad}$ a  convex set, $\varphi(u):=\beta \|u\|_{L^1(\D)}$ a convex function and the objective functional coercive. In particular, the optimal solution $\overline{u}$ satisfies the variational inequality (\cite[Proposition 2.2]{ekeland1999convex})
\begin{equation}\label{eq:variational_inequality}
(\nu \overline{u} -S^\star(y_d-S(\overline{u}+f)),\overline{u} -v)+\varphi(\overline{u})-\varphi(v)\geq 0,\quad \forall v\in U_{ad}.
\end{equation}

Through a pointwise discussion of the box constraints and an analysis of a Lagrange multiplier belonging to the subdifferential of $\varphi$ in $\overline{u}$, \cite{stadler2009elliptic} showed that \eqref{eq:variational_inequality} can be equivalently formulated as the nonlinear equation $\mathcal{F}(\overline{u})=0$, with $\mathcal{F}:L^2(\D)\rightarrow L^2(\D)$ defined as
\begin{equation}\label{eq:nonlinear_optimality_condition}
\medmuskip=-1mu
\thinmuskip=-1mu
\thickmuskip=-1mu
\nulldelimiterspace=0.9pt
\scriptspace=0.9pt 
\arraycolsep0.9em 
\mathcal{F}(u):=u-\frac{1}{\nu}\left(\max(0,\mathcal{T}u-\beta)+\min(0,\mathcal{T}u+\beta)-\max(0,\mathcal{T}u-\beta-\nu b)-\min(0,\mathcal{T}u +\beta-\nu a)\right),
\end{equation}
where $\mathcal{T}:L^2(\D)\ni u \rightarrow -S^\star(Su)+ S^\star(y_d-Sf)\in L^2(\D)$.
Notice that $\mathcal{F}$ is nonsmooth due to the presence of the Lipschitz functions $\max(\cdot)$ and $\min(\cdot)$. 
Nevertheless, $\mathcal{F}$ can be shown to be semismooth \cite{hinze2008optimization}, provided that $\mathcal{T}$ is continuously Fr\'{e}chet differentiable, and further Lipschitz continuous interpreted as map from $L^2(\D)$ to  $L^r(\D)$, with $r>2$ \cite{doi:10.1137/1.9781611970692,hinze2008optimization}. These conditions are satisfied also in our settings since $\mathcal{T}$ is affine and further the adjoint variable $p_\omega$, solution of \eqref{eq:adjoint_equation} with $z=y_d-S(u+f)$, lies in $L^2(\Omega,H^1_0(\D))$ so that $\mathcal{T}u=\E\LQ p_\omega\RQ\in H^1_0(\D)\subset L^r(\D)$, where $r>2$ follows from Sobolev embeddings.

Hence, to solve \eqref{eq:nonlinear_optimality_condition} we use the semismooth Newton method whose iteration reads for $k=1,2,\dots$ until convergence,
\begin{equation}\label{eq:semismooth_problem}
u^{k+1}=u^{k}+du^k,\quad \text{with}\quad \mathcal{G}(u^k)du^k=-\mathcal{F}(u^k),
\end{equation}
$\mathcal{G}(u):L^2(\D)\rightarrow L^2(\D)$ being the generalized derivative of $\mathcal{F}$.
Using the linearity of $\mathcal{T}$ and considering the supports of the weak derivatives of $\max(0,x)$ and $\min(0,x)$, we obtain that
\begin{equation}\label{eq:generalized_differential}
\mathcal{G}(u)[v]=v+\frac{1}{\nu}\chi_{(I^+\cup I^-)}S^\star Sv,
\end{equation}
where $\chi$ is the charateristic function of the union of the disjoint sets 
\[I^+=\left\{x\in \D: 0\leq \mathcal{T}u-\beta\leq \nu b \right\}\text{ and } I^-=\left\{x\in \D: \nu a\leq \mathcal{T}u+\beta\leq 0\right\}.\] 

It is possible to show that the generalized derivative $\mathcal{G}(u)$ is invertible with bounded inverse for all $u$, the proof being identical to the deterministic case treated in \cite{Stadler2}. This further implies that the semismooth Newton method \eqref{eq:semismooth_problem} converges locally superlinearly \cite{doi:10.1137/1.9781611970692}.
We briefly summarize these results in the following proposition. 
\begin{proposition}\label{eq:superlinear_convergence}
Let the initialization $u^0$ be sufficiently close to the solution $\overline{u}$ of \eqref{eq:l1_OCP}. Then the iterates $u^k$ generated by \eqref{eq:semismooth_problem} converge superlinearly to $\overline{u}\in L^2(\D)$.
\end{proposition}

Introducing the supporting variables $dy^k_\omega$ and $dp^k_w$ in $L^2(\Omega;H^1_0(\D))$, the semismooth Newton equation $\mathcal{G}(u^k)du^k=-\mathcal{F}(u^k)$ may be rewritten as the equivalent saddle point system
\begin{equation}\label{full_space_optimality_systeMemismooth1}
\begin{aligned}
& a_\omega(dy^k_\omega,v)-(du^k,v)=0,\quad \forall v\in V,\quad \PP\text{-a-e. } \omega\in \Omega,\\
& a_\omega(v,dp^k_\omega)+(dy^k_\omega,v)=0,\quad \forall v \in V,\ \PP\text{-a-e. } \omega\in \Omega,\\
& (\nu\ du^k - \chi_{(I^+\cup I^-)}\E\LQ dp^k_\omega\RQ,v)_{L^2(\D)}=-\mathcal{F}(u^k),\quad \forall v \in L^2(\D).
\end{aligned}
\end{equation}
Further, if we set $y^0=S(f+u^0)$ and $p^0=S^\star(y_d-y^0)$, due to the linearity of $S$ and $S^\star$, it holds $y^{k+1}=S(u^{k+1})=y^k+dy^{k}$ and similarly $p^{k+1}=p^k+dp^k$.
Once fully discretized and using the notation $\EAP\LQ p_\omega\RQ=\sum_{j=1}^N \zeta_j \pb_{j}$, the optimality condition \eqref{eq:nonlinear_optimality_condition} can be expressed through the nonlinear finite-dimensional map $\Fb:\setR^{N_h}\rightarrow \setR^{N_h}$,
\begin{equation}
\begin{aligned}
\Fb(\ub)&=\ub-\frac{1}{\nu}\Bigl(\max(0,\EAP\LQ \pb_\omega\RQ-\beta)+\min(0,\EAP\LQ \pb_\omega\RQ+\beta)\\
&-\max(0,\EAP\LQ \pb_\omega\RQ-\beta-\nu b)-\min(0,\EAP\LQ \pb_\omega\RQ +\beta-\nu a)\Bigl),
\end{aligned}
\end{equation}
where the $\max(\cdot)$ and $\min(\cdot)$ functions act componentwise.
Equation \eqref{full_space_optimality_systeMemismooth1} leads to the saddle point system
\begin{equation}\label{eq:full_space_optimality_systeMemismooth2}
\begin{pmatrix}
M & & & & A_1^\top\\
& \ddots & & & &\ddots\\
& & M & & & & A_N^\top\\
& & & M & -\zeta_1  M H^k&\dots & -\zeta_N  M H^k\\
A_1 & & & -M\\
& \ddots & &\vdots\\
& & A_N & -M
\end{pmatrix}
\begin{pmatrix}
\dyb^k_1\\ \vdots \\\dyb^k_N\\ \mathbf{du}^k\\
\dpb^k_1\\ \vdots \\ \dpb^k_N
\end{pmatrix}=
\begin{pmatrix}
\mathbf{0}\\ \vdots \\
\mathbf{0}\\ -\Fb(\ub^{k}) \\ \mathbf{0}\\ \vdots \\\mathbf{0}
\end{pmatrix},
\end{equation}
where $H^k\in \mathbb{R}^{N_h\times N_h}$ is a diagonal matrix representing the charateristic function $\chi_{I_k^+\cup I_k^-}$, namely 
\[(H^k)_{i,i}=\frac{1}{\nu} \text{ if }i\in I_k^+\cup I_k^-\quad \text{ and }\quad (H^k)_{i,i}=0 \text{ if }i\notin I_k^+\cup I_k^-,\]
with 
\begin{equation}\label{eq:activesets}
I_k^+=\left\{i: 0\leq \EAP\LQ \pb^k\RQ -\beta\leq \nu b \right\}\text{ and }  I_k^-=\left\{i: \nu a\leq \EAP\LQ \pb^k\RQ +\beta\leq 0\right\}.
\end{equation}
To derive the expression of $H$, we assumed that a Lagrangian basis is used for the finite element space.
Notice that \eqref{eq:full_space_optimality_systeMemismooth2} fits into the general form \eqref{eq:matrix_saddle_point}, and thus we use the collective multigrid algorithm to solve it. Further, with the notation of \eqref{eq:matrix_saddle_point}, it holds
\[
(G)_{i,i}+d_i^\top \diag{a_i}^{-1}\diag{c_i}\diag{a_i}^{-1}e_i=
(M)_{i,i}+(M)^3_{i,i} \sum_{j=1}^N \zeta_j  (A_j)^{-2}_{i,i}>0 \]
if $i\in I^+\cup I^-$, and  
\[(G)_{i,i}+d_i^\top \diag{a_i}^{-1}\diag{c_i}\diag{a_i}^{-1}e_i=(M)_{i,i}>0,\]
if $i\notin I^+\cup I^-$.
The collective multigrid iteration is then well-defined.

The overall semismooth Newton Algorithm is summarized in Algorithm \ref{Alg:semismoothnewton}. At each iteration we solve \eqref{eq:full_space_optimality_systeMemismooth2} using the collective multigrid algorithm  (line 4) and update the active sets given the new iteration (line 10).
Notice that in order to globalize the convergence, we consider a line-search step (lines 6-8) performed on the merit function $\phi(\ub)=\sqrt{\Fb(\ub)^\top M \Fb(\ub)}$ \cite{martinez1995inexact}.
\begin{algorithm}
\setlength{\columnwidth}{\linewidth}
\caption{Globalized semismooth Newton Algorithm to solve $\Fb(\ub)=0$}\label{Alg:semismoothnewton}
\begin{algorithmic}[1]
\Require $\mathbf{u}^0$, $\text{Tol}\in \mathbb{R}^+$, $\sigma,\rho\in (0,1)$. 
\State $\mathbf{y}^0_j=A_j^{-1}(M(\mathbf{f}+\mathbf{u}^0))$, $\mathbf{p}^0_j=\left(A^\top_j\right)^{-1}(M(\mathbf{y}_d-\mathbf{y}^0_j))$, $j=1,\dots,N$.\\
Set $k=0$ and define $I_0^+$ and $I_0^-$ using \eqref{eq:activesets}. 
\While {$\phi(\ub^k)>\text{Tol}$}
\State Solve \eqref{eq:full_space_optimality_systeMemismooth2} calling Alg. \ref{Alg:algorithm-Vcycle} until convergence.
\State Set $\gamma=1$
\While {$\phi(\ub^k+\gamma\dub^k)-\phi(\ub^k)>-\sigma\phi(\ub^k)$}
\State $\gamma=\rho\gamma$.
\EndWhile
\State Update $\mathbf{u}^{k+1}=\mathbf{u}^{k}+\gamma \mathbf{du}^k$, $\mathbf{y}^{k+1}_j=\mathbf{y}_j^{k}+\gamma\mathbf{dy}_j^k$, $\mathbf{p}^{k+1}_j=\mathbf{p}_j^{k}+\gamma\mathbf{dp}_j^k$, $j=1,\dots,N$.
\State Update $I_k^+$ and $I_k^-$ using \eqref{eq:activesets}.
\State Set $k=k+1$.
\EndWhile\\
\Return $\mathbf{u}^{k}, \mathbf{y}_j^{k}$ and $\mathbf{p}_j^{k}$, $j=1,\dots,N$.
\end{algorithmic}
\end{algorithm}

\subsection{Numerical experiments}
In this section we test the semismooth Newton algorithm for the solution of \eqref{eq:nonlinear_optimality_condition} and the collective multigrid algorithm to solve the related optimality system \eqref{eq:full_space_optimality_systeMemismooth2}.
We consider the random PDE-constraint \eqref{eq:state_equation2} with the random diffusion coefficient \eqref{eq:log_normal} set on the L-squared domain. The semismooth iteration is stopped when $\phi(\ub^k)<10^{-9}$. The inner linear solvers are stopped when the relative (unpreconditioned) residual is smaller than $10^{-11}$.

Table \ref{Tab:l1} reports the number of semismooth Newton iterations and in brackets the averaged number of iterations of the V-cycle algorithm used as a solver (left) or as preconditioner for GMRES (right).
Table \ref{Tab:l1} confirms the effectiveness of the multigrid algorithm, which requires essentially the same computational effort as in the linear-quadratic case. 

\begin{table}\caption{Number of semismooth Newton iterations (left), and average number of V-cycle (center) and preconditioned GMRES (right) iterations (in brackets).}\label{Tab:l1}
\centering
\begin{tabular}{| c | c | c | c | c|}
\hline
$\sigma^2$ &  0.1 & 0.5 &  1 & 1.5\\ \hline\hline
It. & 4 \; $\vert$ \;22.5\; $\vert$ \;14 & 5 \; $\vert$ \;22.6\; $\vert$ \;14.2 & 8\; $\vert$ \;23.0\; $\vert$ \;11.8 & 14.9\; $\vert$ \;22.9\; $\vert$ \;15.0 \\
\hline
\end{tabular}\\
\centering
$N_h=705$, $\nu=10^{-4}$ $\beta=10^{-2}$, $N=125$, $N_L=3$, , $L^2=0.5$, $b=50$, $a=-50$.\\\vspace*{0.2cm}
\centering
\begin{tabular}{| c | c | c | c |}
\hline
$N_h$($N_L$) & 161 (2) &  705 (3) & 2945 (4)\\ \hline\hline
It. & 5\; $\vert$ \; 22.0\; $\vert$ \;15.2 &  5\; $\vert$ \;22.6\; $\vert$ \;14.2 & 5 \; $\vert$ \;22.2\; $\vert$ \;14.0 \\
\hline
\end{tabular}\\
\centering
 $\nu=10^{-4}$, $\beta=10^{-2}$, $N=125$, $\sigma^2=0.5$, $L^2=0.5$, $b=50$, $a=-50$.\\\vspace{0.2cm}
\centering
\begin{tabular}{| c | c | c | c | c|}
\hline
$N$ & 8 & 27 & 64 & 125\\ \hline\hline
It. & 5\; $\vert$ \; 21.0\; $\vert$ \;13.0 &  5\; $\vert$ \; 21.6\; $\vert$ \;14.0 &  5\; $\vert$ \; 22.0\; $\vert$ \;14.0 & 5\; $\vert$ \;22.6\; $\vert$ \;14.2 \\
\hline
\end{tabular}\\
\centering
$N_h=705$, $\nu=10^{-4}$, $\beta=10^{-2}$, $\sigma^2=0.5$, $L^2=0.5$, $b=50$, $a=-50$.\\\vspace{0.2cm}
\centering
\begin{tabular}{| c | c | c | c | c|}
\hline
$\beta$ & 0 & $10^{-4}$ & $10^{-3}$ & $10^{-2}$\\ \hline\hline
It. & 4\; $\vert$ \;22.5\; $\vert$ \;14.8 &  4\; $\vert$ \;22.5\; $\vert$ \;14.5 &  5\; $\vert$ \;22.4\; $\vert$ \;14.8 & 5\; $\vert$ \;22.6\; $\vert$ \;14.2 \\
\hline
\end{tabular}\\
\centering
$N_h=705$, $\nu=10^{-4}$, $N=125$, $\sigma^2=0.5$, $L^2=0.5$, $b=50$, $a=-50$.\\\vspace{0.3cm}
\end{table}

More challenging is the limit $\nu\rightarrow 0$ reported in Table \ref{Tab:l1_nu}. The performance of both the (globalized) semismooth Newton iteration and the inner multigrid solver deteriorates. The convergence of the outer nonlinear algorithm can be improved by performing a continuation method, namely we consider a sequence of $\nu=10^{-j}$, $j=2,\dots,8$, and we start the $j$-th problem using as initial condition the optimal solution computed for $\nu=10^{-j+1}$.
Concerning the inner solver, the stand-alone multigrid algorithm struggles since for small values of $\nu$ the optimal control is of bang-bang type, that is satisfies $u=a$, $u=b$ or $u=0$ for almost every point of the mesh (for $\nu=10^{-8}$ only five nodes are nonactive at the optimum). The matrices $H^{k}$ are then close to zero, and the multigrid hierarchy struggles to capture changes at such small scale.
Nevertheless, the multigrid algorithm remains a very efficient preconditioner for GMRES even in this challenging limit.

\begin{table}\caption{Number of semismooth Newton iterations, of V-cycle iterations and of preconditioned GMRES iterations (in brackets). In the second row, the semismooth Newton method starts from a warm-up initial guess obtained through continuation.}\label{Tab:l1_nu}
\centering
\begin{tabular}{| c | c | c | c | c|}
\hline
$\nu$ &  $ 10^{-2}$ & $10^{-4}$ &  $10^{-6}$ & $10^{-8}$\\ \hline\hline
It. & 2\; $\vert$ \;23.0\; $\vert$ \;14.5 &  5\; $\vert$ \;22.7\; $\vert$ \; 14.2 &  17\; $\vert$ \; 25.6\; $\vert$ \; 15.0 &  50\; $\vert$ \; 41.4\; $\vert$ \; 17.2  \\ 
It. & 2\; $\vert$ \;23.0\; $\vert$ \; 14.5 & 4\; $\vert$ \; 22.7\; $\vert$ \; 14.2 & 5\; $\vert$ \;22.25\; $\vert$ \;15.4 & 8\; $\vert$ \; 58.8\; $\vert$ \; 20.9\\
\hline
\end{tabular}\\
\centering
$N_h=705$, $N=125$, $N_L=3$, $\sigma^2=0.5$, $L^2=0.5$, $\beta=10^{-2}$, $b=50$, $a=-50$.\\\vspace*{0.2cm}
\end{table}

\begin{figure}
\includegraphics[scale=0.25]{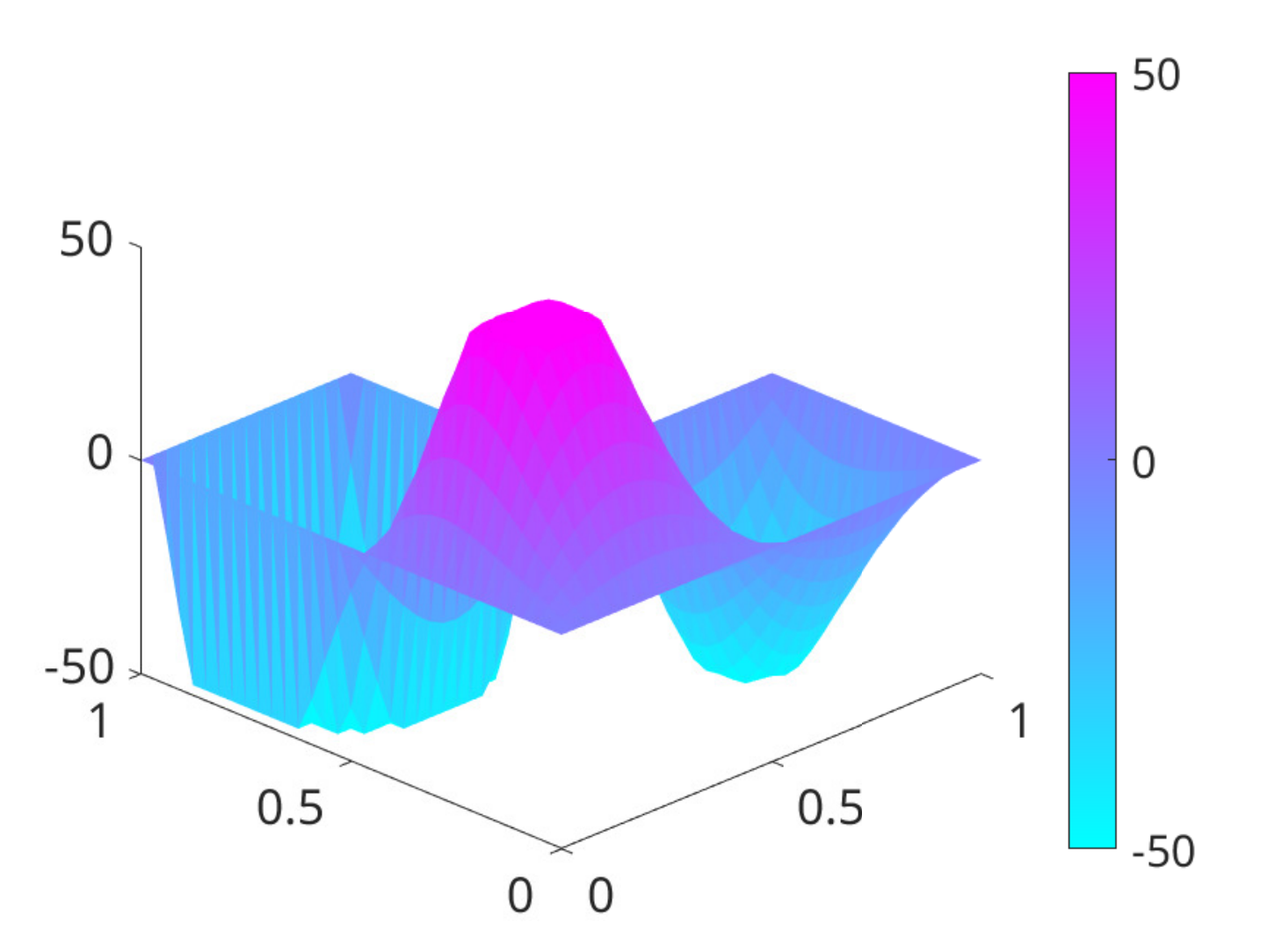}
\includegraphics[scale=0.25]{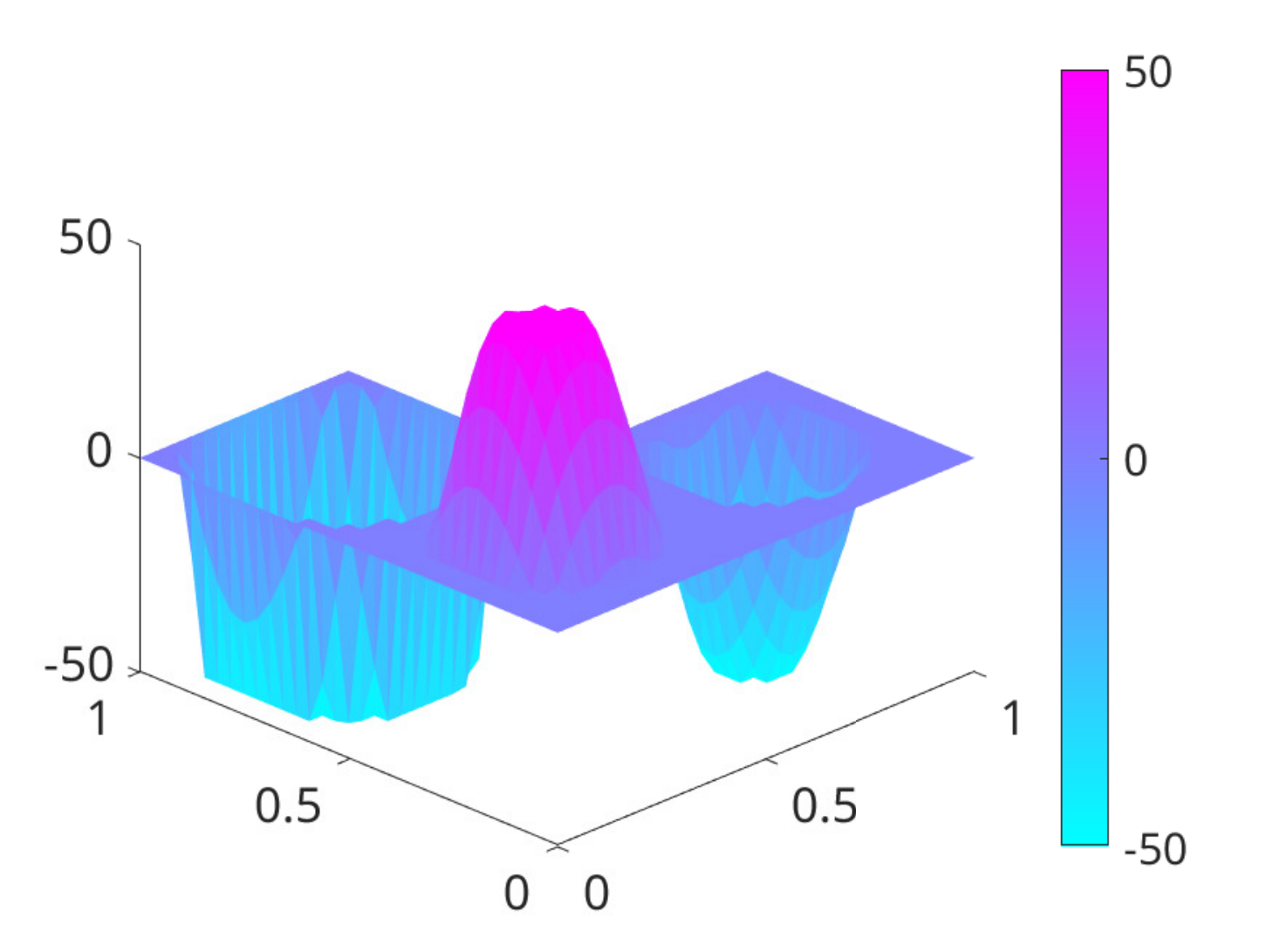}
\includegraphics[scale=0.25]{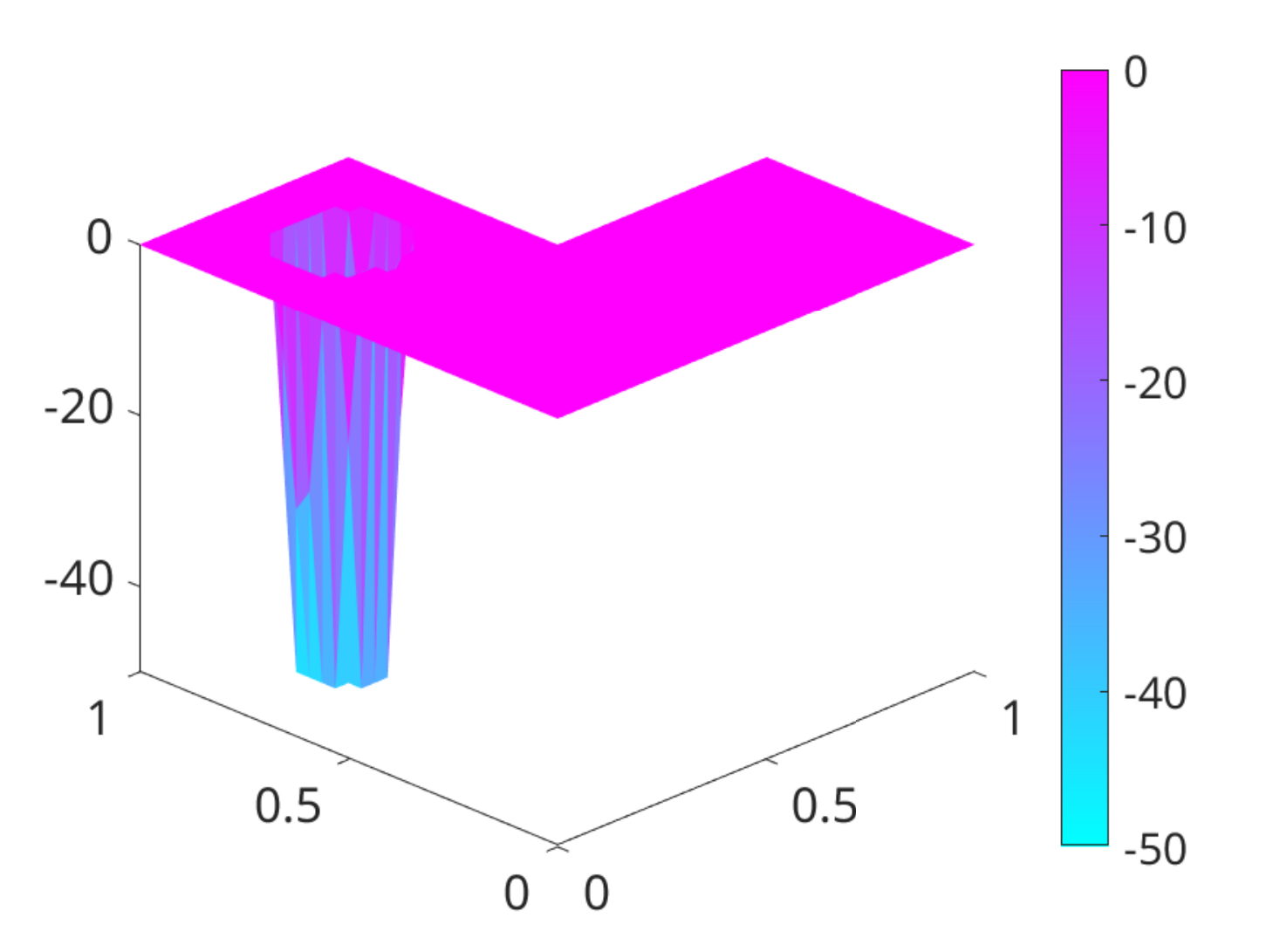}\\
\includegraphics[scale=0.25]{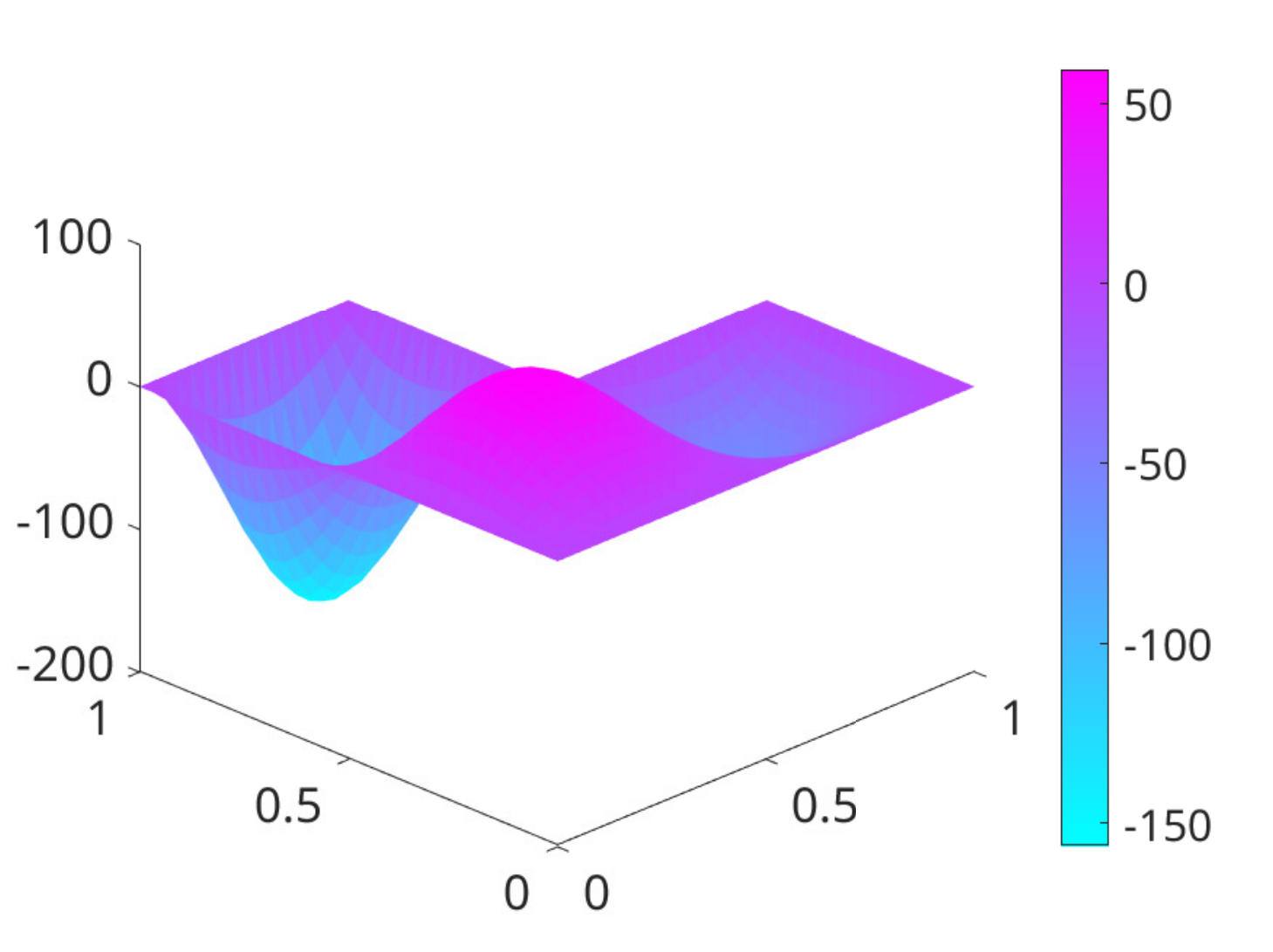}
\includegraphics[scale=0.25]{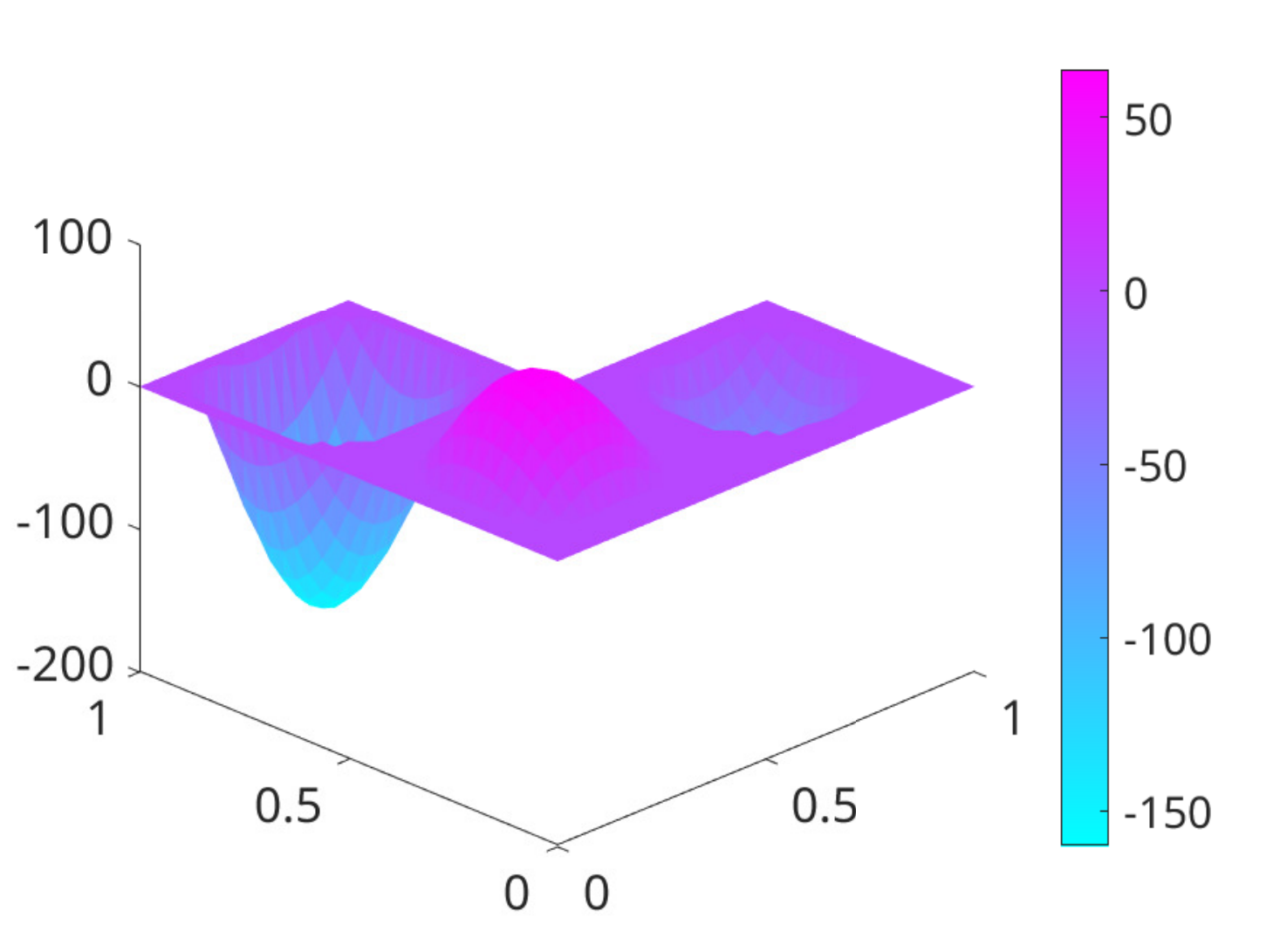}
\includegraphics[scale=0.25]{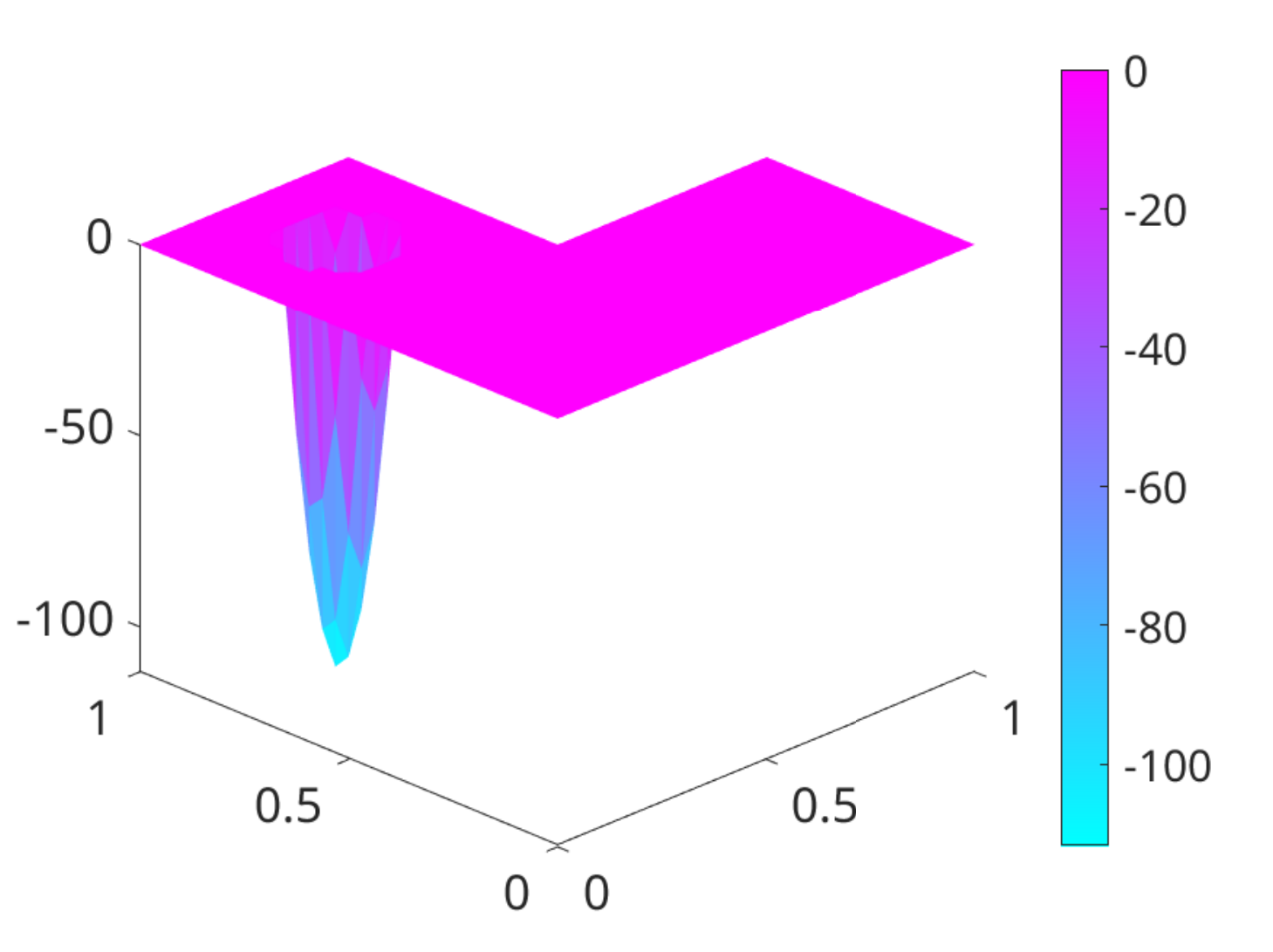}
\caption{From left to right: optimal control computed for $\beta\in \left\{0,5\cdot 10^{-3},5\cdot 10^{-2}\right\}$ with (top row) and without (bottom row) box constraints: $a=-50$, $b=50$.}\label{Fig:betas}
\end{figure}

Fig. \ref{Fig:betas} shows a sequence of optimal controls for different values of $\beta$ with and without box-constraints. The optimal control for $\beta=0$ and without box-constraints corresponds to the minimizer of the linear-quadratic OCP \eqref{eq:quadratic_OCP}. We observe that $L^1$ penalization indeed induces sparsity, since the optimal controls are more and more localized as $\beta$ increases.
Numerically we have verified that for sufficiently large $\beta$, the optimal control is identically equal to zero, a property shown in \cite{stadler2009elliptic}.

\section{A risk-averse optimal control problem under uncertainty}\label{Sec:risk}
In this section we consider an instance of risk-averse OCPUU. This class of problems has recently drawn lot of attention since in engineering applications it is important to compute a control that minimizes the quantity of interest even in rare, but often troublesome, scenarios \cite{Kouri_Cvar,Kouri_ex,antil2021ttrisk,kouri2022primal}.
As a risk-measure \cite{shapiro2014lectures}, we use the Conditional Value-At-Risk (CVaR) of confidence level $\lambda \in (0,1)$,
\begin{equation*}
 \Cvar{X}:= \E\LQ X \vert X\geq \Var{X}\RQ,\quad \forall X\in L^1(\Omega;\setR),
\end{equation*}
that is, the expected value of a quantity of interest $X$ given that the latter is greater than or equal to its $\lambda$-quantile, here denoted by $\Var{X}$.
Rockafellar and Uryasev \cite{rockafellar2000optimization} proved that $\Cvar{X}$ admits the equivalent formulation 
\begin{equation}\label{eq:Cvarequidef}
\Cvar{X}= \inf_{t\in \setR}\left\{t+\frac{1}{1-\lambda}\E\LQ (X-t)^+\RQ \right\},
\end{equation}
where $(\cdot)^+:=\max(0,\cdot)$, if the distribution of $X$ does not have an atom at $\Var{X}$. In order to use tools from smooth optimization, we 
rely on a smoothing approach proposed in \cite{Kouri_Cvar}, which consists in replacing $(\cdot)^+$ with  a smooth function $g_\varepsilon$, $\varepsilon\in\setR^+$, such that $g_\varepsilon\rightarrow (\cdot)^+$ in some functional norm as $\varepsilon\rightarrow 0$.
Specifically, we choose the $C^2$-differentiable approximation
\begin{equation}\label{eq:smoothapproximation}
g_\varepsilon(x)=\begin{cases}
0\quad &\text{if } x\leq -\frac{\varepsilon}{2},\\
\frac{(x-\frac{3}{2})^3}{\varepsilon^2}-\frac{(x-\frac{\epsilon}{2})^4}{2\varepsilon^3}\quad &\text{if } x\in(-\frac{\varepsilon}{2},\frac{\varepsilon}{2}),\\
x\quad &\text{if }x\geq \frac{\varepsilon}{2}.
\end{cases}
\end{equation}
Then, the smoothed risk-averse OCPUU is
\begin{equation}\label{eq:OCP}
\begin{aligned}
&\min_{u\in L^2(\D),t\in \setR} t+\frac{1}{1-\lambda} \E\LQ g_\varepsilon\left(\frac{1}{2}\|y_\omega-y_d\|^2_{L^2(\D)}-t)\right)\RQ+\frac{\nu}{2}\|u\|^2_{L^2(\D)},\\
&\quad \text{subject to}\\
&a_\omega(y_\omega,v)=(u+f,v)\quad \forall v\in V,\ \PP\text{-a.e. }\omega \in \Omega,
\end{aligned}
\end{equation}
where $\nu \in \setR^+$ and $\lambda\in [0,1)$.
The well-posedness of \eqref{eq:OCP}, the differentiability of its objective functional, as well as bounds for the error introduced by replacing $(\cdot)^+$ with $g_{\varepsilon}(\cdot)$, have been analyzed in \cite{Kouri_Cvar}. Further, defining $Q_\omega=\frac{1}{2}\|y_\omega-y_d\|^2_{L^2(\D)}-t$, the optimality conditions form the nonlinear system,
\begin{equation}\label{eq:optmality_systeMmoothed}
\begin{array}{r l r l}
& a_\omega(v,p_\omega)-\frac{g^\prime_\varepsilon\left(Q_\omega\right)}{1-\lambda}(y_d-y_\omega,v)=0,\quad &\forall v \in V,\ \PP\text{-a-e. } \omega\in \Omega,\\
&(\nu\ u -\E\LQ p_\omega\RQ,v)=0,\quad &\forall v\in L^2(\D),\\
&a_\omega(y_\omega,v)-(u+f,v)=0,\quad &\forall v\in V,\quad \PP\text{-a-e. } \omega\in \Omega,\\
&1-\frac{1}{1-\lambda}\E\LQ g^\prime_\varepsilon\left(Q_\omega)\right)\RQ=0.\quad &
\end{array}
\end{equation}

Approximating $V$ and $\E$ with $V_h$ and $\EAP$, and letting $\xbt=(\yb,\ub,\pb,t)$, the finite-dimensional discretization of \eqref{eq:optmality_systeMmoothed} correponds to the nonlinear system $\Fbt(\xbt)=\mathbf{0}$, where $\Fbt:\mathbb{R}^{(2N+1)N_h+1}\rightarrow \mathbb{R}^{(2N+1)N_h+1}$,
\begin{equation}\label{eq:discretized_nonlinearsystem}
\Fbt(\xbt)=\begin{pmatrix}
\Fbt_1(\xbt)\\\Fbt_2(\xbt)\\\Fbt_3(\xbt)\\ \widetilde{F}_4(\xbt)
\end{pmatrix}=\begin{pmatrix}
\widetilde{M}(\yb-I\yb_d)+A^\top \pb\\ \nu M \ub - M\EAP\LQ \pb\RQ\\ A\yb- M(I\ub+\mathbf{f})\\ 1-\frac{1}{1-\lambda}\EAP\LQ g^\prime_\varepsilon(Q_\omega)\RQ
\end{pmatrix},
\end{equation}
with $A=\diag{A_1,\dots,A_N}$, $I=[I_{Nh},\dots,I_{Nh}]\in\mathbb{R}^{N_h\times N_h N}$, $I_h$ being the identity matrix, $\mathbf{y}_d$ is the discretization of $y_d$, and \[\widetilde{M}=\text{diag}\left(\frac{g^\prime_\varepsilon(Q_{\omega_1})}{1-\lambda}M,\dots,\frac{g^\prime_\varepsilon(Q_{\omega_N})}{1-\lambda}M\right), \text{ with }Q_{\omega_j}:=\frac{1}{2}(\mathbf{y}_j-\mathbf{y}_d)^\top M(\mathbf{y}_j-\mathbf{y}_d)-t,\]
for $j=1,\dots,N$. 

A possible approach to solve \eqref{eq:discretized_nonlinearsystem} is to use a Newton method, which given $\xb^k=(\yb^k,\ub^k,\pb^k,t^k)$ computes the corrections $\dxbt^k=(\dyb^k,\dub^k,\dpb^k,dt^k)$ solution of $\Jbt^k\dxbt^k=-\Fbt(\xbt^k)$, where
\begin{equation}\label{eq:optmality_systeMmoothed_discretized}
\Jbt^k:=\begin{pmatrix}
C_1(\mathbf{y}_1^k,t^k) & & & & A_1^\top & & &-\mathbf{v}^k_1\\
& \ddots & & & &\ddots & &\vdots\\
& & C_N(\mathbf{y}_N^k,t^k) & & & & A_N^\top & -\mathbf{v}^k_N\\
& & & \nu M & -\zeta_1 M &\dots & -\zeta_N  M & \\
A_1 & & & -M & & & & \\
& \ddots & &\vdots & & & &  \\
& & A_N & -M & & & &\\
-\zeta_1\left(\mathbf{v}_1^k\right)^\top & \ddots & -\zeta_N\left(\mathbf{v}^k_N\right)^\top & & & & & \frac{\EAP\LQ g_\varepsilon^{\prime\prime}(Q_\omega^k)\RQ}{1-\lambda}
\end{pmatrix},
\end{equation}
with 
\begin{align}
Q_{\omega_i}^k&:=\frac{1}{2}(\mathbf{y}_i^k-\mathbf{y}_d)^\top M(\mathbf{y}_i^k-\mathbf{y}_d) -t^k,\\
C_i(\mathbf{y}_i^k,t^k)&:=\frac{1}{1-\lambda}\left(g^\prime_\varepsilon(Q_{\omega_i}^k)M +g_\varepsilon^{\prime\prime}(Q_{\omega_i}^k)M(\mathbf{y}_i^k-\mathbf{y}_d)(\mathbf{y}_i^k-\mathbf{y}_d)^\top M\right),\label{eq:equation_C}\\
\mathbf{v}^k_i&:=\frac{1}{1-\lambda}g_\varepsilon^{\prime\prime}(Q_{\omega_i}^k)M(\mathbf{y}_i^k-\mathbf{y}_d),
\end{align}
for $i=1,\dots, N$.
Unfortunately, $\Jbt^k$ can be singular away from the optimum, in particular whenever $\EAP\LQ g_\varepsilon^{\prime\prime}(Q_\omega^k)\RQ=0$ which implies
\begin{align}\label{eq:condition}
g_\varepsilon^{\prime\prime}\left(\frac{1}{2}(\mathbf{y}_j^k-\mathbf{y}_d)^\top M(\mathbf{y}_j^k-\mathbf{y}_d) -t^k\right)=0,\ \forall j=1,\dots,N,
\end{align}
which is not unlikely for small $\varepsilon$ since $\text{supp}(g_\varepsilon^{\prime\prime})=(-\frac{\varepsilon}{2},\frac{\varepsilon}{2})$. Splitting strategies have been proposed (e.g. \cite{Markowski2022} in a reduced approach), in which whenever \eqref{eq:condition} is satisfied, an intermediate value of $t$ is computed by solving $\widetilde{F}_4(t;\yb,\ub,\pb)=0$ so to violate \eqref{eq:condition}. 
In the next section, we discuss a similar splitting approach. To speed up the convergence of the outer nonlinear algorithm, we use a preconditioned Newton method based on nonlinear elimination \cite{doi:10.1137/S106482759325154X}. At each iteration we will need to invert saddle-point matrices like \eqref{eq:matrix_saddle_point}, possibly several times. To do so, we rely on the collective multigrid algorithm. 

\subsection{Nonlinear preconditioned Newton method}
Nonlinear elimination is a nonlinear preconditioning technique based on the identification of variables and equations of $\Fb$ (e.g. strong nonlinearities) that slow down the convergence of Newton method. These components are then eliminated through the solution of a local nonlinear problem at every step of an outer Newton. This elimination step provides a better initial guess for the outer iteration, so that a faster convergence is achieved \cite{doi:10.1137/S106482759325154X,doi:10.1137/15M104075X}.

In light of the possible singularity of $\Jbt$, we split the discretized variables $\xbt$ into $\xbt=(\xb,t)$, and we aim to eliminate the variables $\xb$ to obtain a scalar nonlinear equation only for $t$. To do so, we partition \eqref{eq:optmality_systeMmoothed} as
\begin{equation}\label{eq:partitioned}
\Fbt\begin{pmatrix}
\xb\\ t
\end{pmatrix}=\begin{pmatrix}
\Fb_1(\xb,t)\\ F_2(\xb,t)
\end{pmatrix}=\begin{pmatrix}
\mathbf{0}\\0
\end{pmatrix},
\end{equation}
where $\Fb_1=(\Fbt_1(\xb,t),\Fbt_2(\xb,t),\Fbt_3(\xb,t))$ and $F_2(\xb,t)=\widetilde{F}_4(\xb,t)$.
Similarly, $\Jbt$ is partitioned into 
\[\Jbt=\begin{pmatrix}
\Jb_{1,1} & \Jb_{1,2}\\
\Jb_{2,1} & J_{2,2}
\end{pmatrix}\]
whose blocks have dimensions $\Jb_{1,1}\in \setR^{(2N+1)N_h\times (2N+1)N_h}$, $\Jb_{1,2}\in \setR^{(2N+1)N_h\times 1}$, $\Jb_{2,1}\in \setR^{1\times (2N+1)N_h}$, and $\Jb_{2,2}\in \setR$. 
Notice that $\Jb_{1,1}$ is always nonsingular, while $\Jb_{2,1}$, $\Jb_{1,2}$ and $J_{2,2}$ are identically zero if \eqref{eq:condition} is verified.
Thus $\Fb_1$ allows us to define an implicit map $h:\setR\rightarrow \setR^{(2N+1)N_h}$, such that $\Fb_1(h(t),t)=0$, so that the first set of nonlinear equations in \eqref{eq:partitioned} are satisfied. We are then left to solve the nonlinear scalar equation
\begin{equation}\label{eq:reduced_nonlinear}
F(t)=0,\quad\text{where}\quad F(t):=F_2(h(t),t).
\end{equation}
To do so using the Newton method, we need the derivative of $F(t)$ evaluated at $t=t^k$ which, using implicit differentiation, can be computed as
\[F^\prime(t^k)=J_{2,2}(h(t^k),t^k)-\Jb_{2,1}(h(t^k),t^k)\left(\Jb_{1,1}(h(t^k),t^k)
\right)^{-1}\Jb_{1,2}(h(t^k),t^k).\]

The nonlinear preconditioned Newton method is described in Alg. \ref{Alg:preconditioned_newton}, and consists in solving \eqref{eq:reduced_nonlinear} with Newton method. However, to overcome the possible singularity of $J^k_{2,2}$, $\Jb^k_{1,2}$ and $\Jb^k_{2,1}$, we check at each iteration $k$ if \eqref{eq:condition} is satisfied, and in the affirmative case we update $\xb^k$ by solving $\Fb_1(\xb^{k+1},t^k)=0$ using Newton method, and update $t^k$ by solving $F_2(\xb^k,t^{k+1})=0$.
Notice further, that each iteration of the backtracking line-search requires to solve $F_1(h(t),t)=0$ using Newton method, thus additional linear systems with matrix $\Jb_{1,1}$ must be solved.

We report that we also tried to eliminate $t$ by computing the map $l$ such that $F_2(\xb,l(\xb))=0$, while iterating on the variable $\xb$. This has the advantage that $l$ can be evaluted very cheaply, being a scalar equation. However, we needed many more iterations both of the outer Newton method, and consequently of the inner linear solver. Thus, according to our experience, this second approach was less efficient and appealing.

\begin{algorithm}
\setlength{\columnwidth}{\linewidth}
\caption{Nonlinear preconditioned Newton method to solve $\Fbt(\xbt)=0$.}\label{Alg:preconditioned_newton}
\begin{algorithmic}[1]
\Require $t^0$, $\text{Tol}\in \mathbb{R}^+$, $\sigma,\rho\in (0,1)$. 
\State Compute $\xb^0=h(t^0)$ solving $\Fb_1(\xb^0;t^0)=0$ using the Newton method.
\State Set $k=0$.
\While {$\vert F(t^k)\vert >\text{Tol}$}
\If {\eqref{eq:condition} is satisfied}
\State Compute $\xb^{k+1}$ and $t^{k+1}$ solving $\Fb_1(\xb^{k+1};t^k)=0$ and $F_2(\xb^{k+1};t^{k+1})=0$.
\Else
\State Compute Newton's direction $d=-(F^\prime(t^k))^{-1}F(t^k)$.\\
\hspace{1.1cm}Set $\gamma=1$ and compute $\xb=h(t^k+\gamma d)$ solving $\Fb_1(\xb;t^k+d)=0$.
\While {$\vert F(t^k+\gamma d)\vert-\vert F(t^k)\vert >-\sigma\vert F(t^k)\vert $}
\State Set $\gamma=\rho\gamma$.
\State Compute $\xb=h(t^k+\gamma d)$ solving $\Fb_1(\xb;t^k+\gamma d)=0$.
\EndWhile
\State Set $t^{k+1}=t^k+\gamma d$, $\xb^{k+1}=\xb$, $k=k+1$.
\EndIf
\EndWhile\\
\Return $t^{k+1}$ and $\xb^{k+1}$.
\end{algorithmic}
\end{algorithm}

\subsection{Numerical experiments}
In this section we report numerical tests to asses the performance of the preconditioned Newton algorithm to solve \eqref{eq:reduced_nonlinear}, and of the collective multigrid algorithm to invert the matrix $\Jb_{1,1}$.
We consider the random PDE-constraint \eqref{eq:state_equation2} with the random diffusion coefficient \eqref{eq:log_normal}.
Table \ref{Tab:cvar} reports the number of outer and inner Newton iterations, and the average number of V-cycle iterations and of preconditioned GMRES iterations to solve the linear systems at each (inner/outer) Newton iterations. The outer Newton iteration is stopped when $\vert F(t^k)\vert \leq 10^{-6}$, the inner Newton method to compute $h(\cdot)$ is stopped when $\max\left(\|\Fb_{1,1}(\xb^k;t)\|_2/\|\Fb_{1,1}(\xb^0;t)\|_2,\|\Fb_{1,1}(\xb^k;t)\|_2\right)\leq 10^{-8}$, and the linear solvers are stopped when the relative (unpreconditioned) residual is smaller than $10^{-9}$.

In Table \ref{Tab:cvar}, the number of outer Newton iterations is stable, while the number of inner Newton iterations varies between five and fifteen iterations per outer iteration. This is essentially due to how difficult it is to compute the nonlinear map $h(t)$ by solving $\Fb_1(\xb;t)=0$ in line (5), (8) and (11) of Alg. \ref{Alg:preconditioned_newton}. The average number of inner linear solver iterations is quite stable across all experiments. The most challenging case is the limit $\varepsilon\rightarrow 0$ in which we used the solution to the optimization problem as a warmed-up initial guess for the next smaller value of $\varepsilon$. Further, we emphasize that the top left blocks of $\Jb_{1,1}$ involve the matrices $C_i(\yb^k_i,t^k)$ (see \eqref{eq:equation_C}) which contain a dense low-rank term if $g_\varepsilon^{\prime\prime}(Q^k_{\omega_i})\neq 0$. As $\varepsilon\rightarrow 0$, $g_\varepsilon^{\prime\prime}(\cdot)$ tends to a Dirac delta, so the dense term become dominant.
Multigrid methods based on a pointwise relaxations are expected to be not very efficient for these matrices which may not be diagonally dominant.
The standard V-cycle algorithm indeed suffers, however the Krylov acceleration performs better as it handles these low-rank perturbation with smaller effort. For $\varepsilon=10^{-4}$, we sometimes noticed that the GMRES residual stagnates after 20/30 iterations around $10^{-7}/10^{-8}$, due to a loss of orthogonality in the Krylov subspace, and thus resulting in higher number of iterations. We allowed a maximum number of 80 iterations per linear system.

\begin{figure}[h]
\centering
\includegraphics[scale=0.33]{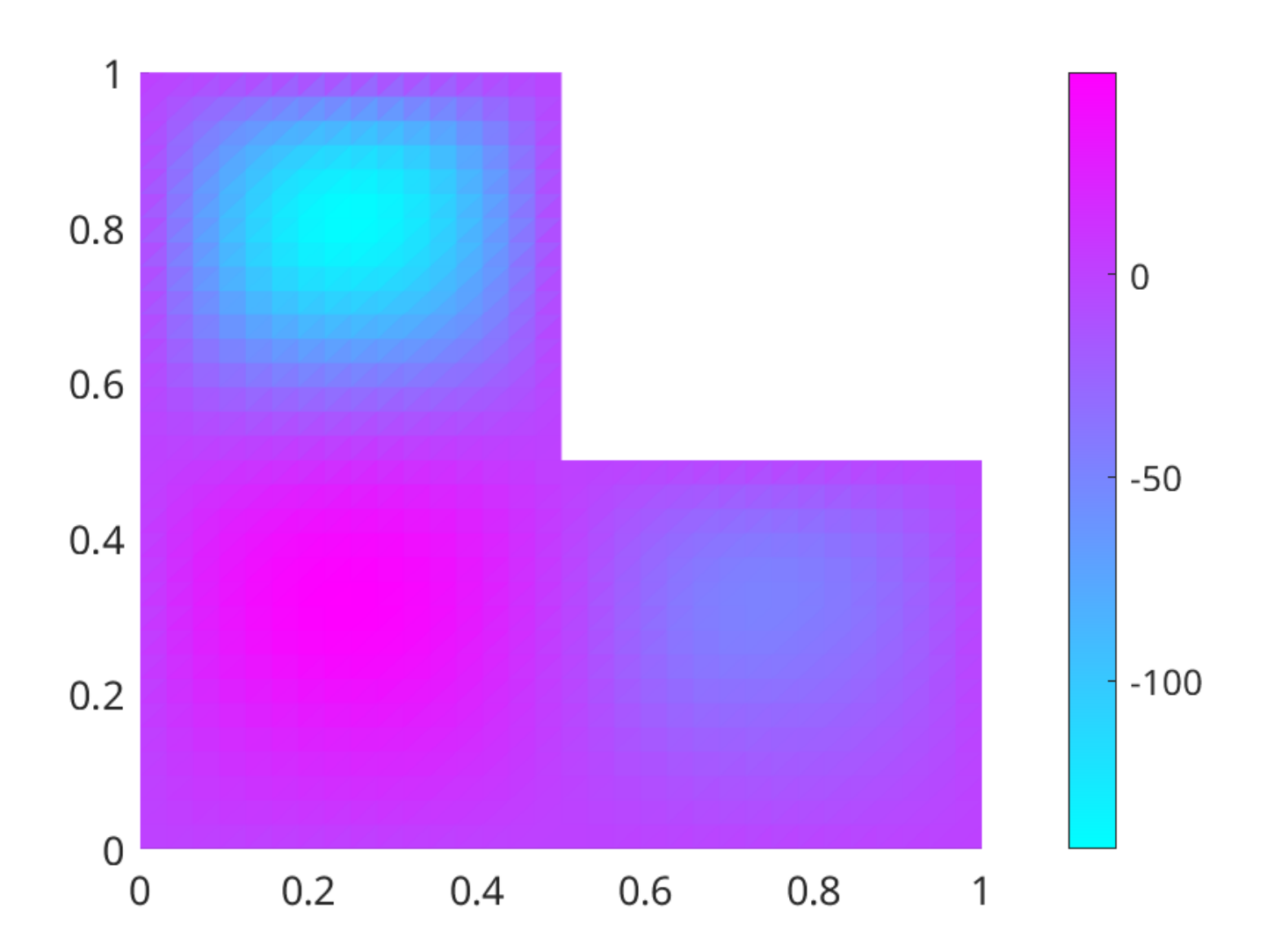}
\includegraphics[scale=0.33]{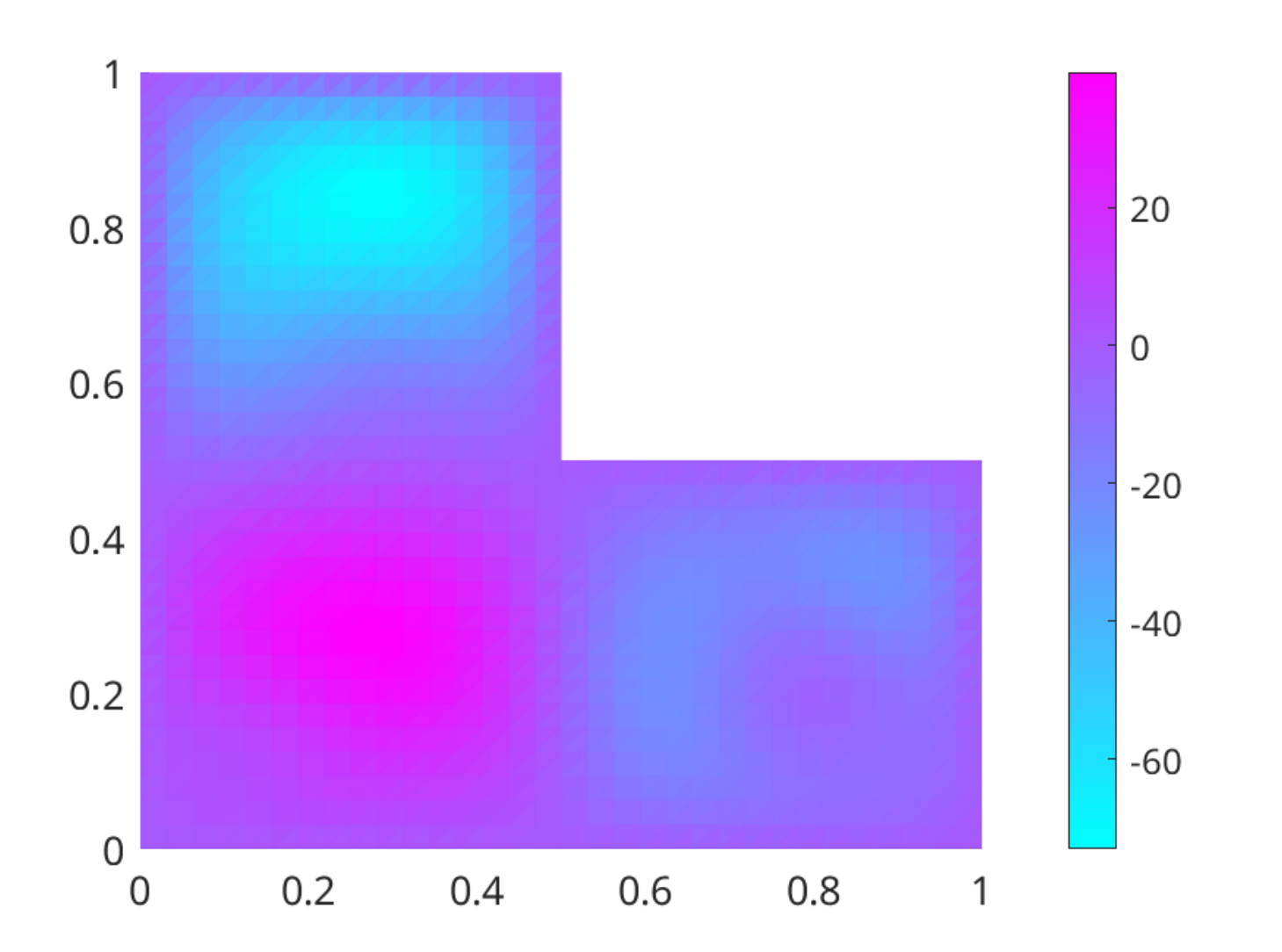}
\includegraphics[scale=0.41]{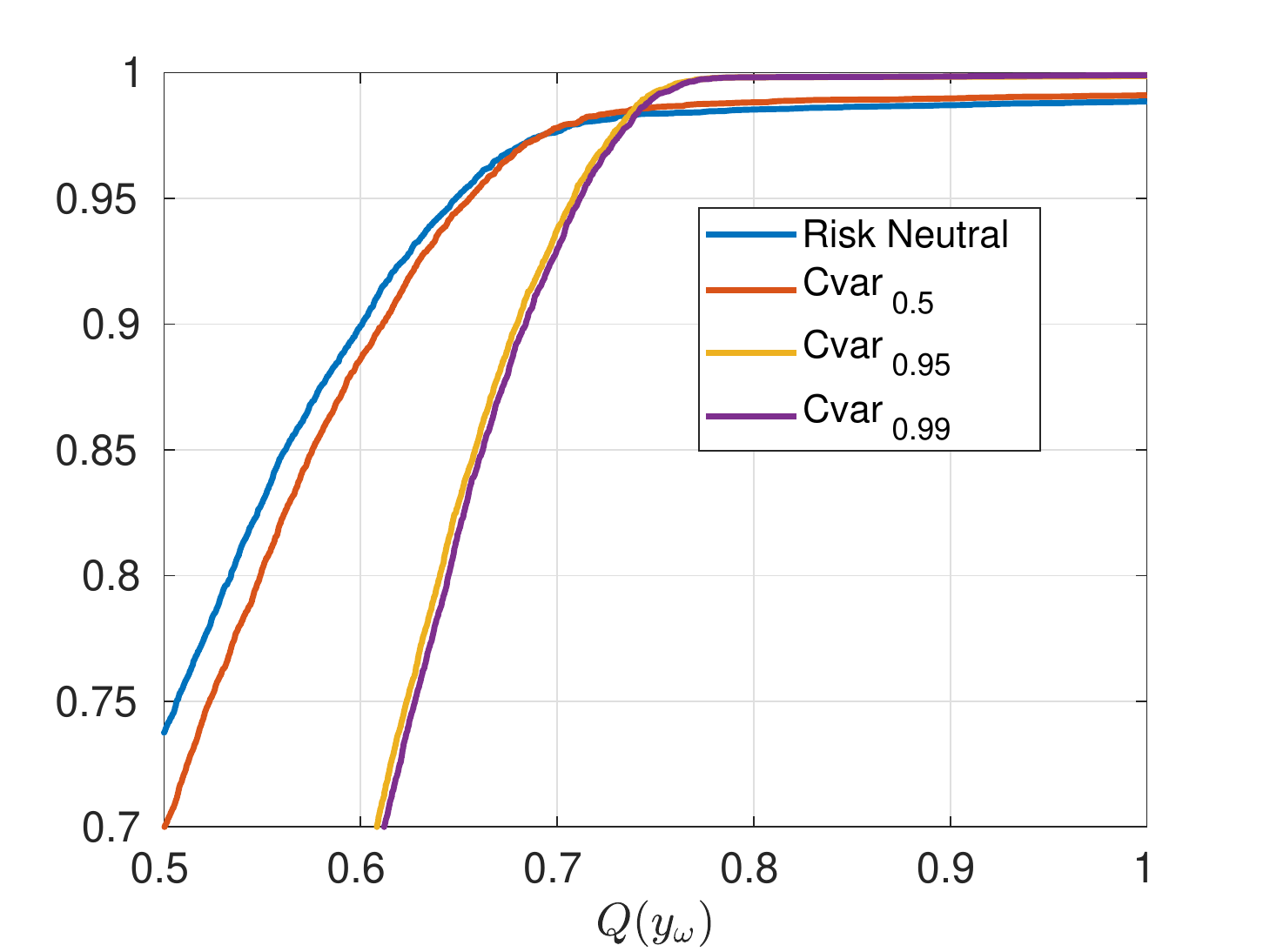}
\caption{Solution of the linear-quadratic OCP (top-left), solution of the smoothed risk-averse OCP with $\lambda=0.99$ (top-right), and cumulative distribution function of the quantity of interest for the controls computed with $\lambda\in \left\{0,0.5,0.95,0.99\right\}$.}\label{Fig:qoi}
\end{figure}

\begin{table}\caption{For each numerical experiment, we report from the left to the right: the number of outer preconditioned Newton iterations, the total number of inner Newton iterations, the averaged number of V-cycle iterations and the averaged number of preconditioned GMRES iterations.}\label{Tab:cvar}
\centering
{\footnotesize
\centering 
\begin{tabular}{| c | c | c | c |}
\hline
$N_h$($N_L$) & 161 (2) &  705 (3) & 2945 (4)\\ \hline\hline
It. & 5\; $\vert$ \;62\; $\vert$ \;23.0\; $\vert$ \;13.9 &  6\; $\vert$ \;79\; $\vert$ \;28.0\; $\vert$ \;15.5 & 6\; $\vert$ \;79\; $\vert$ \;26.2\; $\vert$ \;14.8 \\
\hline
\end{tabular}\\
\centering
 $\nu=10^{-4}$, $N=500$, $\lambda=0.9$, $\varepsilon=10^{-2}$,  $\sigma^2=1$, $L^2=0.1$.\\\vspace{0.2cm}
\centering
\begin{tabular}{| c | c | c | c |}
\hline
$N$ & 500 & 1000 & 2000\\ \hline\hline
It. & 6\; $\vert$ \;63\; $\vert$ \;55.4\; $\vert$ \;17.5 & 5\; $\vert$ \;66\; $\vert$ \;24.4\; $\vert$ \;14.0 &  4\; $\vert$ \;51\; $\vert$ \;24.4\; $\vert$ \;14.0 \\
\hline
\end{tabular}\\
\centering
$N_h=705$, $\nu=10^{-4}$, $\lambda=0.95$, $\varepsilon=10^{-2}$, $\sigma^2=1$, $L^2=0.1$.\\\vspace{0.2cm}
\centering
\begin{tabular}{| c | c | c | c | c|}
\hline
$\lambda$ & 0 & 0.5 & 0.95 & 0.99\\ \hline\hline
It. & 0\; $\vert$ \;1\; $\vert$ \;21.0\; $\vert$ \;14.0 &  5\; $\vert$ \;21\; $\vert$ \;19.4\; $\vert$ \;13.6 &  5\; $\vert$ \;64\; $\vert$ \;23.2\; $\vert$ \;13.8 & 8\; $\vert$ \;129\; $\vert$ \;33.4\; $\vert$ \;17.5 \\
\hline
\end{tabular}\\
\centering
$N_h=705$, $N=2000$, $\nu=10^{-4}$, $\varepsilon=10^{-2}$, $\sigma^2=1$, $L^2=0.1$.\\\vspace{0.3cm}
\centering
\begin{tabular}{| c | c | c | c | c|}
\hline
$\varepsilon$ & $10^{-1}$ & $10^{-2}$ & $10^{-3}$ & $10^{-4}$\\ \hline\hline
It. & 7\; $\vert$ \;67\; $\vert$ \;22.5\; $\vert$ \;17.0 & 3\; $\vert$ \;42\; $\vert$ \;29.1\; $\vert$ \;14.8  &  2\; $\vert$ \;20\; $\vert$ \;$>80$\; $\vert$ \;27.9 & 1\; $\vert$ \;15\; $\vert$ \;58.0 \; $\vert$ \; 55.6  \\
\hline
\end{tabular}\\
\centering
$N_h=705$, $N=1000$, $\nu=10^{-4}$, $\beta=0.95$, $\sigma^2=1$, $L^2=0.1$.\\\vspace{0.3cm}
}
\end{table}
Figure \ref{Fig:qoi} compares the two optimal controls obtained minimizing either $\E\LQ Q(y_\omega)\RQ$ or $\text{CVaR}_{0.99}\LQ Q(y_\omega)\RQ$, and the cumulative distribution functions of $Q(y_{\omega_j})$ computed on 8000 out-of-sample realizations. The risk-averse control indeed minimizes the \textit{risk} of having large values of $Q(y_\omega)$.
The CVaR of level $\lambda=0.99$ is respectively $\text{CVaR}_{0.99}\left( Q(y_{\omega})\right)=2.79$ for the risk-neutral control and $\text{CVaR}_{0.99}\left( Q(y_{\omega})\right)=0.90$ for the risk-averse control.

\section{Conclusion}\label{Sec:conc}
We have presented a multigrid method to solve the large saddle point linear systems that typically arise in full-space approaches to solve OCPUU. We further derived a detailed convergence that fully characterizes the spectrum of the two-level iteration matrix.
The algorithm has been tested as an iterative solver and as a preconditioner on three test cases: a linear-quadratic OCPUU, a nonsmooth OCPUU, and a risk-averse nonlinear OCPUU. Overall, the multigrid method shows very good performances and robustness with respect to the several parameters of the problems considered.

\section{Declarations}
\textbf{Ethical Approval:} Not applicable.

\noindent \textbf{Data Availability:} The codes used in this study are available from the corresponding author on reasonable request.

\noindent \textbf{Competing interests:} The authors declare that they have no conflict of interest.

\noindent \textbf{Funding:} Not applicable.

\noindent \textbf{Authors' contributions:} All authors contributed equally to the manuscript.

\noindent \textbf{Acknowledgements:}The authors wish to thank an anonymous reviewer for the recommendation to develop a convergence analysis of the multigrid algorithm.
G. C. and T. V. are members of GNCS (Gruppo Nazionale
per il Calcolo Scientifico) of INdAM.
The present research is part of the activities of “Dipartimento di Eccellenza 2023-2027".

\end{document}